\documentclass{siamonline220329}
\usepackage{xparse}
\usepackage{xcolor}
\usepackage{xspace}
\usepackage[margin=1.5in]{geometry}
\usepackage{parskip}
\usepackage{amssymb,amsfonts,amsmath}
\usepackage[mathscr]{eucal}
\usepackage{amsbsy}
\usepackage{stmaryrd}
\usepackage{bm}
\usepackage{braket}  
\usepackage{scalerel}
\usepackage{mathtools}
\newsiamremark{problem}{Problem}
\usepackage{tcolorbox}
\tcolorboxenvironment{problem}{colback=white,colframe=black,boxrule=0.5pt}
\usepackage{algorithm, algpseudocode}
\DeclareDocumentCommand{\lineref}{m}{\hyperref[#1]{line~\ref*{#1}}}
\DeclareDocumentCommand{\linesref}{mm}{\hyperref[#1]{lines~\ref*{#1}--\ref{#2}}}
\DeclareDocumentCommand{\Lineref}{m}{\hyperref[#1]{Line~\ref*{#1}}}
\DeclareDocumentCommand{\Linesref}{mm}{\hyperref[#1]{Lines~\ref*{#1}--\ref{#2}}}
\usepackage{subcaption}
\usepackage[inline]{enumitem}
\setlist{noitemsep}
\usepackage{tikz}
\usetikzlibrary{topaths}
\usepackage{hyperref}
\hypersetup{
  naturalnames=false,
  hypertexnames=false,
  breaklinks,
  colorlinks = false,
  hidelinks,
}
\usepackage[capitalize,nameinlink]{cleveref}
\definecolor{c1}{HTML}{e31a1c} %
\definecolor{c1a}{HTML}{fb9a99} %
\definecolor{c2}{HTML}{1f78b4} %
\definecolor{c2a}{HTML}{a6cee3} %
\definecolor{c3}{HTML}{33a02c} %
\definecolor{c3a}{HTML}{b2df8a} %
\definecolor{c4}{HTML}{6a3d9a} %
\definecolor{c4a}{HTML}{cab2d6} %
\definecolor{c5}{HTML}{ff7f00} %
\definecolor{c5a}{HTML}{fdbf6f} %
\definecolor{c9}{HTML}{ffff99} %
\definecolor{c7}{HTML}{b15928} %
\definecolor{c8}{HTML}{f781bf} %
\definecolor{c6}{HTML}{999999} %

\renewcommand{\vec}{\operatorname{vec}}

\renewcommand{\equiv}{:=} %

\NewDocumentCommand{\bigO}{g}{\mathcal{O}\IfValueT{#1}{(#1)}}
\NewDocumentCommand{\qtext}{O{\quad}O{#1}m}{#1\text{#3}#2}
\NewDocumentCommand{\Real}{}{\mathbb{R}}
\NewDocumentCommand{\Natural}{}{\mathbb{N}}
\NewDocumentCommand{\FD}{mm}{\frac{\partial #1}{\partial #2}}

\DeclareMathOperator{\reshape}{reshape}
\NewDocumentCommand{\Tr}{s}{\IfBooleanTF{#1}{\vphantom{\intercal}}{\intercal}}
\NewDocumentCommand{\Vc}{O{} m !g t' t"}{%
  \bm{#1{\mathbf{\MakeLowercase{#2}}}}%
  \IfValueT{#3}{_{#3}}%
  \IfBooleanTF{#4}{^{\Tr}}{%
    \IfBooleanT{#5}{^{\Tr*}}}%
}
\NewDocumentCommand{\Fn}{O{} m !g t' t"}{%
  \bm{#1{\mathnormal{\MakeLowercase{#2}}}}%
  \IfValueT{#3}{_{#3}}%
  \IfBooleanTF{#4}{^{\Tr}}{%
    \IfBooleanT{#5}{^{\Tr*}}}%
}
\NewDocumentCommand{\Mx}{O{} m !g t' t"}{
  \bm{#1{\mathbf{\MakeUppercase{#2}}}}%
  \IfValueT{#3}{_{#3}}%
  \IfBooleanTF{#4}{^{\Tr}}{%
    \IfBooleanT{#5}{^{\Tr*}}}%
}
\NewDocumentCommand{\Qm}{O{} m !g t' t"}{%
  \bm{#1{\mathit{\MakeUppercase{#2}}}}%
  \IfValueT{#3}{_{#3}}%
  \IfBooleanTF{#4}{^{\Tr}}{%
    \IfBooleanT{#5}{^{\Tr*}}}%
}
\NewDocumentCommand{\Cm}{O{} m !g t' t"}{%
  \bm{#1{\mathcal{\MakeUppercase{#2}}}}%
  \IfValueT{#3}{_{#3}}%
  \IfBooleanTF{#4}{^{\Tr}}{%
    \IfBooleanT{#5}{^{\Tr*}}}%
}
\NewDocumentCommand{\Tn}{O{} m !g}{%
  \boldsymbol{#1{\mathscr{\MakeUppercase{#2}}}}%
  \IfValueT{#3}{_{#3}}%
}
\NewDocumentCommand{\Qt}{O{} m !g}{%
  \boldsymbol{#1{\mathcal{\MakeUppercase{#2}}}}%
  \IfValueT{#3}{_{#3}}%
}
\NewDocumentCommand{\Tm}{O{} m O{} !g t' t"}{
  \bm{#1{\mathbf{\MakeUppercase{#2}}}}_{\prn#3(#4)}%
  \IfBooleanTF{#5}{^{\Tr}}{%
    \IfBooleanT{#6}{^{\Tr*}}}%
}
\NewDocumentCommand{\megaLRexp}{mmmmmmmm}{%
  \IfBooleanTF{#1}{%
    \IfBooleanTF{#2}{%
      \IfBooleanTF{#3}{%
        \IfBooleanTF{#4}%
        {\Biggl#6 #8 \Biggr#7}%
        {\biggl#6 #8 \biggr#7}}%
      {\Bigl#6 #8 \Bigr#7}}%
    {\bigl#6 #8 \bigr#7}}%
  {\IfBooleanTF{#5}
    {\left#6 #8 \right#7}
    {#6 #8 #7}}%
}
\NewDocumentCommand{\prn}{ ssss t+ g d() t'}
{%
  \IfValueT{#6}{\megaLRexp{#1}{#2}{#3}{#4}{#5}{(}{)}{#6}}%
  \IfValueT{#7}{\megaLRexp{#1}{#2}{#3}{#4}{#5}{(}{)}{#7}}%
  \IfBooleanT{#8}{^{\Tr}}%
}
\NewDocumentCommand{\abs}{ ssss t+ g d||}
{%
  \IfValueT{#6}{\megaLRexp{#1}{#2}{#3}{#4}{#5}{\vert}{\vert}{#6}}%
  \IfValueT{#7}{\megaLRexp{#1}{#2}{#3}{#4}{#5}{\vert}{\vert}{#7}}%
}
\NewDocumentCommand{\crly}{ ssss t+ g }
{%
  \IfValueT{#6}{\megaLRexp{#1}{#2}{#3}{#4}{#5}{\{}{\}}{#6}}%
}
\NewDocumentCommand{\nrm}{ ssss t+ g }
{%
  \IfValueT{#6}{\megaLRexp{#1}{#2}{#3}{#4}{#5}{\|}{\|}{#6}}%
}
\NewDocumentCommand{\fnrm}{ ssss t+ g }
{%
  \IfValueT{#6}{\megaLRexp{#1}{#2}{#3}{#4}{#5}{\|}{\|_{\rm F}}{#6}}%
}
\NewDocumentCommand{\tnrm}{ ssss t+ g }
{%
  \IfValueT{#6}{\megaLRexp{#1}{#2}{#3}{#4}{#5}{\|}{\|_{2}}{#6}}%
}
\NewDocumentCommand{\sqr}{ ssss t+ g d[] }
{%
  \IfValueT{#6}{\megaLRexp{#1}{#2}{#3}{#4}{#5}{\lbrack}{\rbrack}{#6}}%
  \IfValueT{#7}{\megaLRexp{#1}{#2}{#3}{#4}{#5}{\lbrack}{\rbrack}{#7}}%
}
\NewDocumentCommand{\ang}{ ssss t+ g d[] }
{%
  \IfValueT{#6}{\megaLRexp{#1}{#2}{#3}{#4}{#5}{\langle}{\rangle}{#6}}%
  \IfValueT{#7}{\megaLRexp{#1}{#2}{#3}{#4}{#5}{\langle}{\rangle}{#7}}%
}
\NewDocumentCommand{\dsqr}{ ssss t+ g d[] }
{%
  \IfValueT{#6}{\megaLRexp{#1}{#2}{#3}{#4}{#5}{\llbracket}{\rrbracket}{#6}}%
  \IfValueT{#7}{\megaLRexp{#1}{#2}{#3}{#4}{#5}{\llbracket}{\rrbracket}{#7}}%
}

\renewcommand{\top}{\mathbin{\mathchoice
    {\vcenter{\hbox{\scaleobj{0.75}{\varbigcirc}}}}
    {\vcenter{\hbox{\scaleobj{0.75}{\varbigcirc}}}}
    {\vcenter{\hbox{\scaleobj{0.6}{\varbigcirc}}}}
    {\vcenter{\hbox{\scaleobj{0.5}{\varbigcirc}}}}
  }}

\newcommand{\krn}{\mathbin{\otimes}}

\newcommand{\krp}{\mathbin{\odot}}

\newcommand{\had}{\mathbin{\mathchoice
    {\vcenter{\hbox{\scaleobj{1.5}{\ast}}}}
    {\vcenter{\hbox{\scaleobj{1.5}{\ast}}}}
    {\vcenter{\hbox{\scaleobj{1.2}{\ast}}}}
    {\vcenter{\hbox{\scaleobj{1.0}{\ast}}}}
  }}

\NewDocumentCommand{\dsub}{m m t_ m}{#1_{#2_{#4}}}
\NewDocumentCommand{\tttMlist}{m m m m m m m}{%
  #1_{1} #2 %
  \IfBooleanF{#5}{#1_{2} #2} %
  \IfBooleanTF{#6}{ %
    #1_{3} %
    \IfBooleanT{#7}{#2 #1_{4} } %
  }{ %
    \IfBooleanTF{#3}{\cdots}{\dots} #2 #1_{#4} %
  }%
}
\NewDocumentCommand{\tttRMlist}{m m m m m m m}{%
  \IfBooleanTF{#6}%
  {\IfBooleanT{#7}{#1_{4} #2} #1_{3} #2 #1_{2} #2 #1_{1}}%
  {#1_{#4} \IfBooleanF{#5}{#2 #1_{#4-1}} #2 \cdots  #2 #1_{1}}%
}
\NewDocumentCommand{\tttSlist}{m m m m m m}{%
  #1_{1} %
  #2 \IfBooleanTF{#3}{\cdots}{\dots}
  #2 #1_{#5-1} %
  \IfValueT{#6}{#2 #6} %
  #2 #1_{#5+1} %
  #2 \IfBooleanTF{#3}{\cdots}{\dots}
  #2 #1_{#4} %
}
\NewDocumentCommand{\tttRSlist}{m m m m m m}{%
  #1_{#4} %
  #2 \IfBooleanTF{#3}{\cdots}{\dots}
  #2 #1_{#5+1} %
  \IfValueT{#6}{#2 #6} %
  #2 #1_{#5-1} %
  #2 \IfBooleanTF{#3}{\cdots}{\dots}
  #2 #1_{1} %
}
\NewDocumentCommand{\miwc}{s s t! O{i} O{d}}{%
  \tttMlist{#4}{,}{\BooleanFalse}{#5}{#3}{#1}{#2}%
}
\NewDocumentCommand{\siwc}{O{k} O{i} O{d} g}{%
  \tttSlist{#2}{,}{\BooleanFalse}{#3}{#1}{#4}%
}
\NewDocumentCommand{\minc}{s s t! O{i} O{d}}{%
  \tttMlist{#4}{}{\BooleanTrue}{#5}{#3}{#1}{#2}%
}

\NewDocumentCommand{\sinc}{O{k} O{i} O{d} g}{%
  \tttSlist{#2}{}{\BooleanTrue}{#3}{#1}{#4}%
}
\NewDocumentCommand{\tttTsze}{>{\SplitArgument{3}{,}}m}{\tttTszeFour#1} %
\NewDocumentCommand{\tttTszeFour}{m m m m}{#1 \IfValueT{#2}{\times #2} \IfValueT{#3}{\times #3} \IfValueT{#4}{\times #4}} %

\NewDocumentCommand{\msiz}{s s t! O{n} O{d} !g !g}{%
  \IfValueTF{#6}%
  {\tttTsze{#6}\IfValueT{#7}{\times \cdots \times #7}}%
  {#4_1 \IfBooleanF{#3}{\times #4_2}%
    \IfBooleanTF{#1}{\times #4_3 \IfBooleanT{#2}{\times #4_4}}{\times \cdots \times #4_{#5}}}%
}

\NewDocumentCommand{\Rmsiz}{s s t! O{n} O{d} !g !g}{%
  \mathbb{R}^{\IfValueTF{#6}%
  {\tttTsze{#6}\IfValueT{#7}{\times \cdots \times #7}}%
  {#4_1 \IfBooleanF{#3}{\times #4_2}%
    \IfBooleanTF{#1}{\times #4_3 \IfBooleanT{#2}{\times #4_4}}{\times \cdots \times #4_{#5}}}}%
}

\NewDocumentCommand{\ssiz}{O{k} O{n} O{d} g}{%
  \tttSlist{#2}{\times}{\BooleanTrue}{#3}{#1}{#4}%
}

\NewDocumentCommand{\Rssiz}{O{k} O{n} O{d} g}{%
  \mathbb{R}^{\tttSlist{#2}{\times}{\BooleanTrue}{#3}{#1}{#4}}%
}
\NewDocumentCommand{\tttMdom}{>{\SplitArgument{3}{,}}m}{\tttMdomFour#1} %
\NewDocumentCommand{\tttMdomFour}{m m m m}{[#1] \IfValueT{#2}{\otimes [#2]} \IfValueT{#3}{\otimes [#3]} \IfValueT{#4}{\otimes [#4]}} %

\NewDocumentCommand{\mdom}{s s t! O{n} O{d} !g !g}
{
  \IfValueTF{#6}
  {\tttMdom{#6}\IfValueT{#7}{\otimes \cdots \otimes [#7]}}
  {[#4_1] \IfBooleanF{#3}{\otimes [#4_2]}
    \IfBooleanTF{#1}{\otimes [#4_3] \IfBooleanT{#2}{\otimes [#4_4]}}{\otimes \cdots \otimes [#4_{#5}]}}
}

\NewDocumentCommand{\sdom}{O{k} O{n} O{d} g}{%
  [#2_{1}] %
  \otimes \cdots
  \otimes [#2_{#1-1}] %
  \IfValueT{#4}{\otimes [#4]} %
  \otimes [#2_{#1+1}] %
  \otimes \cdots
  \otimes [#2_{#3}] %
}
\NewDocumentCommand{\tttMsummand}{>{\SplitArgument{1}{/}}m}{\tttMsummandTwo#1} %
\NewDocumentCommand{\tttMsummandTwo}{mm}{\IfValueTF{#2}{\sum_{#1=1}^{#2}}{\sum_{#1}}} %

\NewDocumentCommand{\msum}{s s t! O{i} O{n} O{d} >{\SplitList{,}}g !g}{%
  \IfValueTF{#7}%
  {\ProcessList{#7}{\tttMsummand}\IfValueT{#8}{\cdots\tttMsummand{#8}}}%
  {%
    \tttMsummand{#4_1/#5_1}%
    \IfBooleanF{#3}{\tttMsummand{#4_2/#5_2}}%
    \IfBooleanTF{#1}%
    {\tttMsummand{#4_3/#5_3} \IfBooleanT{#2}{\tttMsummand{#4_4/#5_4}}} %
    {\cdots\tttMsummand{#4_{#6}/#5_{#6}}} %
  }%
}
\usepackage{etoolbox}
\usepackage{xcolor}
\usepackage{pgfplots}
\usepgfplotslibrary{groupplots}
\usepackage{pgfplotstable}
\usepgfplotslibrary{colorbrewer}
\usetikzlibrary{external}
\colorlet{c1}{Set1-A}
\colorlet{c2}{Set1-B}
\colorlet{c3}{Set1-C}
\colorlet{c4}{Set1-D}

\tikzset{
  sc/.append style={below,midway,inner ysep=0em,align=flush left},
}

\NewDocumentCommand{\cphifi}{}{{CP-HiFi}\xspace}
\NewDocumentCommand{\KF}{o}{\bm{\mathcal{K}}\IfValueT{#1}{^{(#1)}}}
\NewDocumentCommand{\T}{}{\Tn{T}}
\NewDocumentCommand{\QT}{}{\Qt{T}}
\NewDocumentCommand{\Tobs}{}{\Tn[\bar]{T}}
\NewDocumentCommand{\Rmnp}{}{\Rmsiz{m,n,p}}
\NewDocumentCommand{\va}{}{\Vc{a}}
\NewDocumentCommand{\vb}{}{\Vc{b}}
\NewDocumentCommand{\vc}{}{\Vc{c}}
\NewDocumentCommand{\fa}{}{\Fn{a}}
\NewDocumentCommand{\fb}{}{\Fn{b}}
\NewDocumentCommand{\fc}{}{\Fn{c}}
\NewDocumentCommand{\ff}{}{\Fn{f}}

\NewDocumentCommand{\vx}{}{\Vc{x}}
\NewDocumentCommand{\vy}{}{\Vc{y}}
\NewDocumentCommand{\vv}{}{\Vc{v}}
\NewDocumentCommand{\vw}{}{\Vc{w}}
\NewDocumentCommand{\A}{}{\Mx{A}}
\NewDocumentCommand{\B}{}{\Mx{B}}
\NewDocumentCommand{\C}{}{\Mx{C}}
\NewDocumentCommand{\QA}{}{\Qm{A}}
\NewDocumentCommand{\QB}{}{\Qm{B}}
\NewDocumentCommand{\QC}{}{\Qm{C}}
\NewDocumentCommand{\QK}{}{\Qm[\hat]{K}}
\NewDocumentCommand{\I}{}{\Mx{I}}
\NewDocumentCommand{\K}{}{\Mx{K}}
\NewDocumentCommand{\V}{}{\Mx{V}}
\NewDocumentCommand{\W}{}{\Mx{W}}
\NewDocumentCommand{\Z}{}{\Mx{Z}}
\NewDocumentCommand{\HK}{}{\mathcal{H_{\KF}}}
\NewDocumentCommand{\xk}{O{k}}{x_{#1}}
\NewDocumentCommand{\fakl}{O{k}G{\ell}}{\fa^{(#1)}_{#2}}
\NewDocumentCommand{\vakl}{O{k}G{\ell}}{\va^{(#1)}_{#2}}
\NewDocumentCommand{\vwkl}{O{k}G{\ell}}{\Vc{w}^{(#1)}_{#2}}
\NewDocumentCommand{\spt}{G{k}}{v_{#1}}
\NewDocumentCommand{\sptk}{O{k}G{i_k}}{v^{(#1)}_{#2}}
\NewDocumentCommand{\vk}{O{k}}{\vv^{(#1)}}
\NewDocumentCommand{\Kk}{O{k}}{\Mx{K}_{#1}}

\NewDocumentCommand{\ExpFig}{smm}{
\begin{figure}[htp]
  \tikzexternaldisable
  \centering
  \begin{tikzpicture}
    \IfBooleanTF{
    \node[name=bar,inner sep=0em]
    {\includegraphics{extfigs/quasitensor-2d-exp-#2.pdf}};
    \path let \p1=(bar.south west), \p2=(bar.south east) in (\p1) --
    (\p2) node[sc,text width=0.9*(\x2-\x1),name=foo-caption] {\subcaption{Red dots
        illustrate observed noisy samples per $(i,j)$ pair. Red lines
        show underlying true quasi-tensor. Grey vertical lines show
        $x$-values for which at least one $(i,j)$ pair has a
        sample.}\label{fig:exp-#2-data}};
    }{\node (foo-caption)}
    \path (foo-caption.south) node[inner sep=0em,anchor=north,name=cp]
    {\includegraphics{extfigs/exp-#2-cp-factors.pdf}};
    \path let \p1=(cp.south west), \p2=(cp.south east) in (\p1) --
    (\p2) node[sc,text width=0.9*(\x2-\x1),name=cp-caption] {\subcaption{CP solution}\label{fig:exp-#2-cp}};
    \path (cp-caption.south) node[inner sep=0em,anchor=north,name=hifi]
    {\includegraphics{extfigs/exp-#2-hifi-factors.pdf}};
    \path let \p1=(hifi.south west), \p2=(hifi.south east) in (\p1) --
    (\p2) node[sc,text width=0.9*(\x2-\x1),name=hifi-caption] {\subcaption{CP-HIFI solution}\label{fig:exp-#2-hifi}};
    
  \end{tikzpicture}
  \caption{Experiment \# #2: #3}
  \label{fig:exp-#2}
\end{figure}

}

\begin{document}

\title{Tensor Decomposition Meets RKHS:\\Efficient Algorithms for Smooth and Misaligned Data}
\author{Brett W. Larsen\thanks{Databricks Mosaic Research; Centers for Computational Mathematics and Neuroscience, Flatiron Institute, New York, NY (\email{brettlarsen@flatironinstitue.org}).}%
  \and
  Tamara G.~Kolda\thanks{MathSci.ai, Dublin, CA (\email{tammy.kolda@mathsci.ai}).}%
  \and
  Anru R. Zhang\thanks{Department of Biostatistics and Bioinformatics and Department of Computer
Science, Duke University (\email{anru.zhang@duke.edu})}%
  \and
  Alex H. Williams\thanks{Center for Neural Science, New York University; Center for Computational Neuroscience, Flatiron Institute (alex.h.williams@nyu.edu)}%
}
\maketitle

\begin{abstract}
  The canonical polyadic (CP) tensor decomposition decomposes
  a multidimensional data array into a sum of outer products of
  finite-dimensional vectors.
  Instead, we can replace some or all of the vectors with continuous
  functions (infinite-dimensional vectors) from a reproducing kernel Hilbert space (RKHS).
  We refer to tensors with some infinite-dimensional modes as quasi-tensors,
  and the approach of decomposing a tensor with some continuous RKHS modes is referred to as
  \cphifi (hybrid infinite and finite dimensional) tensor decomposition.
  An advantage of \cphifi is that it can enforce smoothness in the
  infinite dimensional modes. Further, \cphifi does not require the observed data to lie
  on a regular and finite rectangular grid and naturally incorporates misaligned data.
  We detail the methodology and illustrate it on a synthetic example.
\end{abstract}

\section{Introduction}
\label{sec:introduction}

\begin{figure}[ht]
  \centering
  \tikzexternaldisable
  \begin{tikzpicture}
    \node[anchor=north east,name=bar,inner sep=0em]
    {\includegraphics{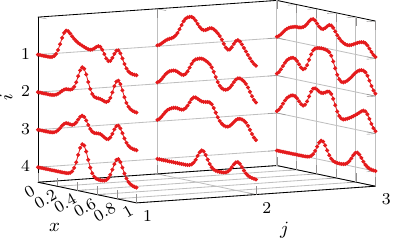}};
    \path let \p1=(bar.south west), \p2=(bar.south east) in
    (\p1) -- (\p2) node[sc,text width=0.9*(\x2-\x1)]
    {\subcaption{Observed data from quasitensor, 50 evenly spaced and aligned samples
        per $(i,j)$ pair.}\label{fig:data}};
    \node[anchor=north west,name=bar2,inner sep=0em]
    {\includegraphics{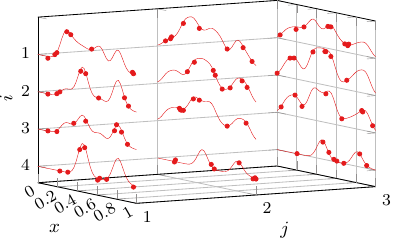}};
    \path let \p1=(bar2.south west), \p2=(bar2.south east) in
    (\p1) -- (\p2) node[sc,text width=0.9*(\x2-\x1)] (bar2-caption)
    {\subcaption{Observed data from quasitensor, 8 unevenly spaced and unaligned samples
        per $(i,j)$ pair.}\label{fig:sampled-data}};
    \begin{scope}[yshift=-6cm]
      \node[anchor=north,name=foo,inner sep=0em]
      {\includegraphics{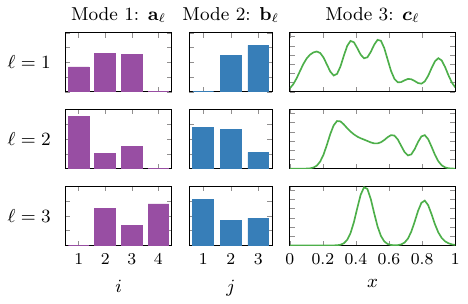}};
    \end{scope}
    \path let \p1=(foo.south west), \p2=(foo.south east) in (\p1) --
    (\p2) node[sc,text width=0.9\textwidth] {\subcaption{Rank-3
        decomposition of quasitensor. Factor matrices corresponding to
        first two modes are finite (vectors) and factor quasi matrix
        corresponding to third mode is infinite
        (functions).}\label{fig:cp-data}};
  \end{tikzpicture}
  \caption{Example $4 \times 3 \times [0,1]$ quasitensor:
    $\QT(i,j,x)$.
  }
  \label{fig:intro-example}
\end{figure}
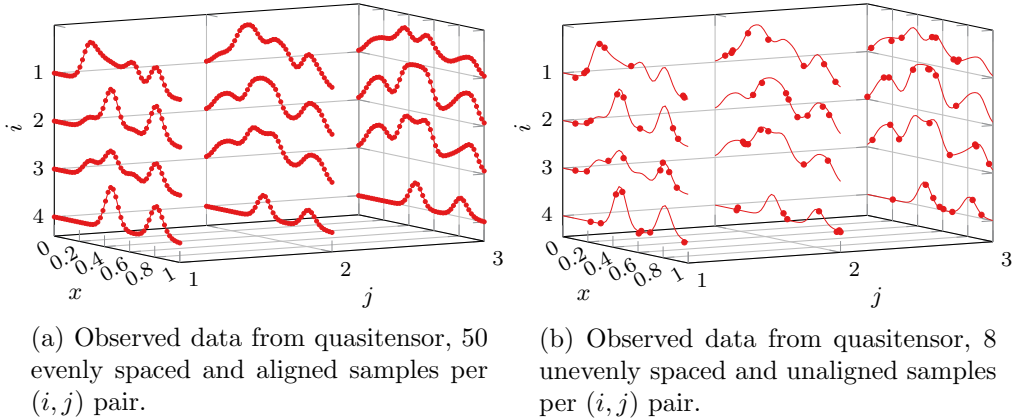
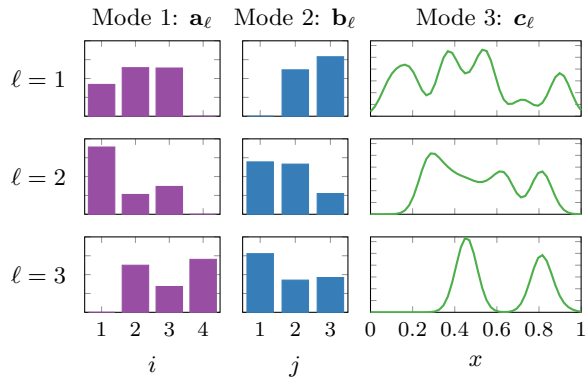

In many real-world situations, a data tensor is the result of making discrete
measurements on a set of continuous processes.
We refer to such a set of continuous processes as a \emph{quasitensor}.
For example, a quasitensor $\QT \in \Rmsiz{m,n,\infty}$ is an $m \times n$
collection of functions.
In \cref{fig:data}, we show an example of a three-way quasitensor,
$\QT \in \Rmsiz{4,3,\infty}$,
measured at discrete regularly spaced samples in the third mode.
In \cref{fig:sampled-data}, we show the same quasitensor measured
at discrete irregularly spaced samples.

The goal of this work is to show how to compute a
CP hybrid infinite and finite dimensional (\cphifi) decomposition.
Considering the example in \cref{fig:intro-example},
we want a rank-$r$ decomposition of the quasitensor $\QT \in \Rmsiz{m,n,\infty}$ as
\begin{equation}\label{eq:hifi-cp-3}
  \tikzexternaldisable
  \begin{tikzpicture}[baseline=(a.base)]
    \path[use as bounding box] node[inner xsep=0em] (a) {$\QT$};
    \draw[gray,<-] (a.north) to[out=90,in=0] ++(120:0.3)
    node[left,inner sep=0.1em,font=\footnotesize]{quasitensor};
  \end{tikzpicture}
  \approx \sum_{\ell=1}^r \,  
  \va{\ell}
  \circ
  \vb{\ell}
  \circ
  \begin{tikzpicture}[baseline=(c.base)]
    \path[use as bounding box] node[inner xsep=0em] (c) {$\fc{\ell}$};
    \draw[gray,<-] (c.north) to[out=90,in=180] ++(60:0.3)
    node[right,inner sep=0.1em,font=\footnotesize]{function};
  \end{tikzpicture}
  \equiv \dsqr{\A,\B,
    \begin{tikzpicture}[baseline=(a.base)]
      \path[use as bounding box] node[inner xsep=0em] (a) {$\QC$};
      \draw[gray,<-] (a.north) to[out=90,in=180] ++(60:0.3)
      node[right,inner sep=0.1em,font=\footnotesize]{quasimatrix};
    \end{tikzpicture}
  }.
\end{equation}
Here, $\A\in\Real^{m \times r}$ and $\B \in \Real^{n \times r}$ are matrices as usual,
with $\va{\ell}$ and $\vb{\ell}$ denoting their $\ell$th columns, respectively.
The difference from the usual CP decomposition
is that $\QC\in \Real^{\infty \times r}$ is a \emph{quasimatrix}, i.e.,
a collection of $r$ functions, and $\fc{\ell}$ denotes the $\ell$th function.
Pointwise, this means
\begin{equation}\label{eq:hcp-3way}
  \QT(i,j,x) \approx \sum_{\ell=1}^r \va{\ell}(i) \, \vb{\ell}(j) \, \Fn{c}{\ell}(x)
  \qtext{for all}
  (i,j,x) \in \mdom{m,n} \otimes \Real.
\end{equation}

For the infinite dimensional functions, we choose 
functions from a reproducing kernel Hilbert space (RKHS).
A positive semidefinite kernel function $\KF: [a,b] \otimes [a,b] \rightarrow \Real$
uniquely defines an RKHS, denoted $\HK$.
Considering the example in \cref{eq:hifi-cp-3}, we would require
each function $\fc{\ell} \in \HK$.
The advantage of using functions from an RKHS is that we can
reduce the minimization over an infinite-dimensional subspace to
an equivalent problem in discrete space via the
representer theorem \cite{ScHeSm01}.

This paper is organized as follows.
In \cref{sec:preliminaries}, we review matrix and tensor concepts, including
quasimatrices and quasitensors.
In \cref{sec:rkhs-repr-theor}, we review reproducing kernel Hilbert spaces
and kernel ridge regression for fitting functions.
Because the general $d$-way tensor notation is cumbersome,
we first introduce \cphifi in terms of 3-way tensors with a single continuous mode
in \cref{sec:simple-3-way}.
This is expanded to the general $d$-way case with any number of continuous modes
in \cref{sec:general-d-way}. The general case includes matrices ($d=2$) as well as
higher-order tensors, and it works with up to $d$ infinite-dimensional modes.
We provide a demonstration of the method in \cref{sec:numerical-examples}.
This problem and variants have been considered in many different scenarios,
and we review related works in \cref{sec:related-work}.
Of particular note, Han, Shi, and Zhang \cite{HaShZh23} and
Tang, Kolda, and Zhang \cite{TaKoZh24} discuss a 3-way model with a single
continuous mode using RKHS and applications to a microbiome dataset. 
Here, we focus on computational and efficiency issues while
also extending the model to general $d$-way tensors.

\section{Preliminaries: Matrices and tensors}
\label{sec:preliminaries}

We denote vectors as boldface lowercase letters,
matrices as boldface upper letters, and tensors
using boldface Euler letters.
We refer to any matrix or tensor with at least one
continuous mode as a quasimatrix or quasitensor, respectively.
We will sometimes distinguish the case where all modes are continuous
by instead using the term cmatrix and ctensor (short for continuous).
We denote functions, quasimatrices, and quasitensors
using slanted variants of the notation for vectors, matrices, and tensors.
Thus, $\va$ is a vector
whereas $\fa$ is a function which we bold to emphasize that this typically requires an infinite dimensional vector to represent. $\A$ is a matrix whereas
$\QA$ is a quasimatrix, and $\T$ is a tensor whereas $\QT$ is a quasitensor.
We let $\Natural = \set{1,2,\dots}$ denote the set of natural numbers.
For any $n \in \Natural$, we use the shorthand $[n]\equiv\set{1,2,\dots,n}$.

\subsection{Matrix operations}
\label{sec:matr-tens-oper}
We let $\tnrm{\cdot}$ denote the vector 2-norm and $\fnrm{\cdot}$ denote the matrix
Frobenious norm.
The vectorization of a matrix $\A \in \Rmsiz{m,n}$ stacks its columns,
i.e., entry $(i,j)$ of $\A$ is mapped to entry $i+m(j-1)$ of $\vec(\A)$.

\paragraph{Weighted norm}

Let $\K \in \Rmsiz{p,p}$ be a positive definite matrix, and let
$\Mx{W}$ be a $p \times n$ matrix. Then we can define
\begin{equation}\label{eq:mat-nrm-K}
  \nrm{\Mx{W}}_{\K}^2 \equiv \sum_{i=1}^n \Vc{w}{i}' \K \Vc{w}{i}"
  = \vec(\W)^{\Tr} (\Mx{I}{n} \krn \K) \vec(\W).
\end{equation}
If $\K$ is the identity matrix, then $\nrm{\Mx{W}}_{\K}^2 = \nrm{\Mx{W}}_F^2$ (the Frobenious norm).

\paragraph{Hadamard product}
Given matrices $\A,\B \in \Rmsiz{m,n}$, their \textbf{Hadamard product} is denoted
$\A \had \B \in \Rmsiz{m,n}$ with $(\A \had \B)_{ij} = a_{ij} b_{ij}$ for all $(i,j)$.

\paragraph{Kronecker product}
Given matrices $\A \in \Rmsiz{m,n}$ and $\B \in \Rmsiz{p,q}$, their
\textbf{Kronecker product} is denoted $\A \krn \B \in \Rmsiz{mn,pq}$ and
defined as
\begin{displaymath}
  \A \krn \B =
  \begin{bmatrix}
    a_{11} \B & \cdots & a_{1n} \B \\
    \vdots & \ddots & \vdots \\
    a_{m1} \B & \cdots & a_{mn} \B \\
  \end{bmatrix}.
\end{displaymath}
\paragraph{Khatri-Rao product}
Given matrices $\A \in \Rmsiz{m,r}$ and $\B \in \Rmsiz{n,r}$
with columns given by $\va{k}$ and $\vb{k}$, their
\textbf{Khatri-Rao product} is denoted $\A \krp \B \in \Rmsiz{mn,r}$ and
defined as
\begin{displaymath}
  \A \krp \B =
  \begin{bmatrix}
    \va{1} \krn \vb{1} & \cdots & \va{r} \krn \vb{r}
  \end{bmatrix}.
\end{displaymath}

The following Kronecker and Khatri-Rao properties are useful for our discussion:
\begin{align}
  \label{eq:krn-transpose}
  \prn(\A \krn \B)' &= \A' \krn \B',\\
  \label{eq:krn-krn}
  (\A \krn \B)(\C \krn \Mx{D}) &= (\A\C) \krn (\B\Mx{D}),\\
  \label{eq:krn-vec}
  \vec(\A\C\B') &= (\B \krn \A) \vec(\C),\\
  \label{eq:krp-had}
  \prn(\A \krp \B)'(\A \krp \B) &= (\A'\A) \had (\B'\B),\\
  \label{eq:krn-krp}
  (\A \krn \B)(\C \krp \Mx{D}) &= (\A\C) \krp (\B\Mx{D}).
\end{align}

\subsection{Matrix calculus}

A few rules of matrix calculus are useful. In the following,
we assume $\A \in \Rmsiz{m,n}$, $\Vc{x} \in \Real^n$, and $\vb \in \Real^m$:
\begin{align}
  \label{eq:mc-1}
  \FD{ \nrm{ \A \Vc{x} - \vb }_2^2 }{\Vc{x}} &= 2\A'(\A\Vc{x}-\vb),\\
  \label{eq:mc-2}
  \FD{ \Vc{x}' \A \Vc{x} }{\Vc{x}} &= (\A+\A')\Vc{x}.
\end{align}

\subsection{Quasimatrices}

A quasimatrix $\QC \in \Real^{\infty \times r}$
represents a collection of $r$ univariate functions defined on the same interval
$\mathcal{I} \in \Real$.
We may write
\begin{displaymath}
  \QC =
  \begin{bmatrix}
    \fc{1} & \fc{2} & \cdots \fc{r}
  \end{bmatrix}
  \qtext{where} \fc{j}: \mathcal{I} \rightarrow \Real
  \qtext[\;]{for all} j \in [r].
\end{displaymath}
Oftentimes, we want to evaluate that set of functions at a set of
points, producing a matrix as the output. This is a way of
going from a quasimatrix to a matrix.
Given a vector $\vv \in \mathcal{I}^p$ (a set of points to be evaluated),
we can define the matrix
\begin{displaymath}
  \QC(\vv) \equiv
  \begin{bmatrix}
    \fc_1(v_1) & \fc_2(v_1) & \cdots & \fc_r(v_1) \\
    \fc_1(v_2) & \fc_2(v_2) & \cdots & \fc_r(v_2) \\
    \vdots & \vdots & \ddots & \vdots \\
    \fc_1(v_p) & \fc_2(v_p) & \cdots & \fc_r(v_p) \\
  \end{bmatrix} \in \Rmsiz{p,r}.
\end{displaymath}

We extend the previous notions of quasimatrix and
consider $\KF \in \Rmsiz{\infty,\infty}$
with two infinite-dimensional modes, representing a bivariate function:
\begin{displaymath}
  \KF: \mathcal{I} \otimes \mathcal{I} \rightarrow \Real.
\end{displaymath}
Such a bivariate function is sometimes called a cmatrix for continuous matrix when considered in the context of quasimatrices.
We can go to a quasimatrix with one infinite-dimensional mode
by evaluating at a set of points in just one of the two modes.
Given a vector $\vv \in \mathcal{I}^p$, we can define the quasimatrix
\begin{displaymath}
  \KF(\cdot,\vv) \equiv
  \begin{bmatrix}
    \KF(\cdot,v_1) & \KF(\cdot, v_2) & \cdots & \KF(\cdot,v_p)
  \end{bmatrix}
  \in \Rmsiz{\infty,p}.
\end{displaymath}
Further, we can evaluate at a set of points in each mode,
going from a quasimatrix to a matrix.
Given $\vv \in \mathcal{I}^p$ and $\Vc{w} \in \mathcal{I}^q$, we can define
the matrix
\begin{displaymath}
  \KF(\vv,\Vc{w}) \equiv
  \begin{bmatrix}
    \KF(v_1,w_1) & \KF(v_1, w_2) & \cdots & \KF(v_1,w_p)\\
    \KF(v_2,w_1) & \KF(v_2, w_2) & \cdots & \KF(v_2,w_p)\\
    \vdots & \vdots & \ddots & \vdots \\
    \KF(v_p,w_1) & \KF(v_p, w_2) & \cdots & \KF(v_p,w_p)\\
  \end{bmatrix}
  \in \Rmsiz{p,q}.
\end{displaymath}
We say $\KF$ is \emph{symmetric} if $\KF(x,y) = \KF(y,x)$ for all $x,y \in \mathcal{I}$.
We say $\KF$ is \emph{positive semidefinite} (psd) if for any $p \in \Natural$ and
$\vv \in \mathcal{I}^p$, the matrix $\KF(\vv,\vv)$ is positive semidefinite.
We note that here that we are using $\KF(\vv,\vv)$ as shorthand for evaluating $\KF$ at all combinations of the entries of $\vv$, not that $\KF$ takes in vector inputs.
\subsection{Tensor operations and tensor decompositions}

Given \emph{factor matrices}
$\A \in \Rmsiz{m,r}$, $\B \in \Rmsiz{n,r}$ and $\C \in \Rmsiz{p,r}$,
we denote their columns as
\begin{displaymath}
  \A = 
  \begin{bmatrix}
    \va{1} & \va{2} & \cdots & \va{r}
  \end{bmatrix}
  , \quad
  \B = 
  \begin{bmatrix}
    \vb{1} & \vb{2} & \cdots & \vb{r}
  \end{bmatrix}
  , \qtext{and}
  \C = 
  \begin{bmatrix}
    \vc{1} & \vc{2} & \cdots & \vc{r}
  \end{bmatrix}
  .
\end{displaymath}
We can express a CP decomposition
using \emph{Kruskal tensor notation} as 
\begin{displaymath}
  \dsqr{\A,\B,\C} = \sum_{\ell=1}^r \va{\ell} \top \vb{\ell} \top \vc{\ell},
\end{displaymath}
where $\top$ denotes the tensor outer product.
Elementwise, this means
\begin{displaymath}
  \prn**{\dsqr{\A,\B,\C}}_{ijk} = \sum_{\ell=1}^r \va_{\ell}(i) \, \vb_{\ell}(j) \, \vc_{\ell}(k)
  \qtext{for all}
  (i,j,k) \in \mdom{m,n,p}.
\end{displaymath}
More generally, for $\A{k} \in \Rmsiz{n_k,r}$ for $k \in[d]$,
the notation $\dsqr{\miwc![\A]}$ means
\begin{displaymath}
  \prn**( \dsqr{\miwc![\A]} )_{\minc!} = \msum{\ell/r} \prod_{k=1}^d \vakl(i_k)
  \qtext{for all}
  (\miwc!) \in \mdom!.
\end{displaymath}
Here $\vakl$ represents column $\ell$ of matrix $\A{k}$.

\paragraph{Vectorization}

The vectorization of a three-way tensor $\T \in \Rmnp$ reorders the entries of the tensor into a
vector  $\vec(\T) \in \Real^{mnp}$ such that entry $(i,j,k)$ of $\T$ is mapped to entry $i + m(j-1) + mn(k-1)$ of $\vec(\T)$.

This can be extended to $d$-way tensors. For $\T \in \Rmsiz$,
we have $\vec(\T) \in \Real^N$ with $N=\prod_{k=1}^d n_k$. For the exact formula for the
mapping, see \cite{KoBa09}.

\paragraph{Mode-$k$ Unfolding}

The mode-$k$ unfolding of a tensor $\T$ is denoted
as $\Tm{T}{k}$ and rearranges the tensor of size $\msiz$ into a
matrix of size $n_k \times N_k$ where $N_k = \prod_{\ell \neq k} n_{\ell}$;
see \cite{KoBa09} for details.
There is a connection between the Khatri-Rao product and Kruskal tensor
unfoldings:
\begin{align}
  \label{eq:ktensor-unfold-1}  
  \prn*{ \dsqr{ \A,\B,\C }}_{(1)} &= \A\prn(\C \krp \B)',\\
  \label{eq:ktensor-unfold-2}
  \prn*{ \dsqr{ \A,\B,\C }}_{(2)} &= \B\prn(\C \krp \A)',\\
  \label{eq:ktensor-unfold-3}
  \prn*{ \dsqr{ \A,\B,\C }}_{(3)} &= \C\prn(\B \krp \A)',\\
  \label{eq:ktensor-unfold}
  \prn*{ \dsqr{ \miwc[\A] }}_{(k)} &= \A{k}\prn(\A{d} \krp \cdots \krp \A{k+1}
  \krp \A{k-1} \krp \cdots \krp \A{1})'.
\end{align}

We define the \emph{mode-$k$ perfect shuffle matrix} to be the matrix $\Mx{P}{k}$ such
that
\begin{equation}\label{eq:tensor-perfect-shuffle}
  \Mx{P}{k} \vec(\T) = \vec\prn{\Tm{T}{k}}.
\end{equation}
So, for example, $\vec(\dsqr{ \A,\B,\C }) = \Mx{P}{2}' \vec(\B\prn(\C \krp \A)')
=\Mx{P}{2}'\sqr*{ (\C \krp \A) \krn \Mx{I}} \vec(\B)$.
The last step used \cref{eq:krn-vec}.

\paragraph{MTTKRP} 

A key kernel for computing CP decomposition is the
\emph{matricized-tensor times Khatri-Rao product (MTTKRP)}
kernel.
This is the product of the mode-$k$ unfolding of a $d$-way tensor with a
Khatri-Rao product of $d-1$ matrices. 
For a three-way tensor $\T \in \Rmsiz{m,n,p}$ and matrices
$\A \in \Rmsiz{m,r}$, $\B \in \Rmsiz{n,r}$, and $\C \in \Rmsiz{p,r}$,
the MTTKRP operations are $\Tm{T}{1} (\C \krp \A)$ (mode 1),
$\Tm{T}{2} (\C \krp \B)$ (mode 2), and $\Tm{T}{3} (\B \krp \A)$ (mode 3).
The computational complexity in any mode is $\bigO{mnpr}$.
For a $d$-way tensor $\T \in \Rmsiz$ and matrices
$\A{k} \in \Rmsiz{n_k,r}$ for $k \in [d]$, the MTTKRP operation
has computational complexity is $\bigO(r N)$ where $N = \prod_{k \in[d]} n_k$
and its form is
\begin{displaymath}
  \Tm{T}{k}\prn(\A{d} \krp \cdots \krp \A{k+1} \krp \A{k-1} \krp \cdots \krp \A{1}).
\end{displaymath}

\paragraph{Norms}

We define the norm of a tensor $\T$ to be
\begin{equation}\label{eq:nrm-def}
  \nrm{\T} \equiv \nrm{\!\vec(\T)}_2 = \nrm{\Tm{T}{k}}_F.
\end{equation}

Let $N$ be the total number of indices in $\T$, 
$\Omega$ be a subset of tensor indices, and $q = \abs{\Omega} < N$ be the
number of indices in $\Omega$.
We define $\Mx{S}{\Omega} \in \Rmsiz{N,q}$ as the operator such that
$\Mx{S}{\Omega}'\vec(\T)$ picks out only the entries of $\T$ in $\Omega$.
The matrix $\Mx{S}{\Omega}$ is an orthonormal matrix consisting of $q$ columns of the $N \times N$ identity matrix.
Finally, we define 
\begin{equation}\label{eq:omega-nrm-def}
  \nrm{\T}_{\Omega}^2 \equiv \sum_{(\miwc) \in \Omega} \T(\miwc)^{2}
  = \nrm{\Mx{S}{\Omega} \vec(\T)}_2^2.
\end{equation}
If $\Omega$ is the set of all indices, then $\nrm{\T}_{\Omega} = \nrm{\T}$.

\section{Preliminaries: Reproducing kernel Hilbert spaces (RKHS)}
\label{sec:rkhs-repr-theor}

As discussed in the introduction, we constrain the functions in our \cphifi decomposition
to live in an RKHS. We review the necessary concepts by
considering the following generic problem.

\begin{problem}[Functional regression] Given a set of \emph{design points}
$\set{x_i}_{i=1}^p \subset \mathcal{I}\equiv[a,b] \subset \Real$ and
corresponding noisy function evaluations $\set{y_i}_{i=1}^p$, find a
function $\ff: \mathcal{I} \rightarrow \Real$ such that
$\ff(x_i) \approx y_i$ for all $i \in [p]$.  
\end{problem}

We must constrain our choices for function $\ff$ in some way.
We could, for instance, constrain $\ff$ to be a linear combination
of a finite number of basis functions.
In this work, we focus on constraining  $\ff$ to lie in an RKHS.
An RKHS is an infinite-dimensional Hilbert space of functions.
An RKHS can be uniquely defined by a
choice of symmetric psd kernel function
$\KF: \mathcal{I} \otimes \mathcal{I} \rightarrow \Real$.
\Cref{tab:kernel-functions} provides a few common examples of kernels defined on $\Real$;
many more psd kernels can be found in \cite[Chapter 3]{FaMc15},
\cite[Chapter 4]{WiRa06}, and \cite{Du14}.
We let $\HK$ denote the RKHS defined by psd kernel $\KF$.
We can formally express the RKHS version of our problem as follows.

\begin{table}[ht]
  \centering\small\tikzexternaldisable
  \renewcommand{\arraystretch}{2}
  \begin{tabular}{|l|c|l|}  
    \hline
    \bf Name & \bf Kernel Function $K(x,y)$ & \bf User Parameters \\ \hline
    Exponential & \tikz[baseline=(a.base)]{\node(a){$\displaystyle\exp\prn+(\frac{-\abs{x-y}}{c})$};} & $c > 0$\\
    Gaussian & \tikz[baseline=(a.base)]{\node(a){$\displaystyle\exp\prn+(\frac{-\prn{x-y}^2}{2c^2})$};} & $c \neq 0$\\
    Rational Quadratic & \tikz[baseline=(a.base)]{\node(a){$\displaystyle\prn+(1+\frac{-\prn{x-y}^2}{2\alpha c^2})^{-\alpha}$};} & $c\neq 0$ and $\alpha>0$\\
    Sinc &  \tikz[baseline=(a.base)]{\node(a){$\displaystyle\frac{\sin\prn+(c(x - z))}{c(x - z)}$};} & $c>0$ \\
    Periodic & \tikz[baseline=(a.base)]{\node(a){$\displaystyle \exp\prn***{ -\frac{2 \sin^2 \prn*( \pi \abs{x-z} / p ) }{2c^2}  }$};} & $p>0$ and $c \neq 0$ \\
    \hline
  \end{tabular}
  \caption{Examples of psd kernel functions}
  \label{tab:kernel-functions}
\end{table}

\begin{problem}[Kernel ridge regression] 
Given a set of \emph{design points} $\set{x_i}_{i=1}^p \subset \mathcal{I}$,
corresponding noisy function evaluations $\set{y_i}_{i=1}^p$,
symmetric psd kernel $\KF: \mathcal{I} \otimes \mathcal{I} \rightarrow \Real$, and
regularization parameter $\lambda>0$, solve
\begin{equation}
  \label{eq:KRR}
  \min_{\ff \in \HK}
  \sum_{i=1}^p  \prn*( y_i - \ff(x_i) )^2 + \lambda \nrm{\ff}_{\HK}^2.
\end{equation}
\end{problem}
Although \cref{eq:KRR} is an optimization problem over an infinite dimensional
space, by the Representer Theorem (see \cref{sec:repr-theor-kern}),
its solution is guaranteed to have the form
\begin{equation}
  \label{eq:KRR-solution-form}
  \ff^*= \sum_{i=1}^p w_i \, \KF(\cdot,x_i) = \QK\vw,
\end{equation}
where $\QK \equiv \KF(\cdot,\Vc{x}) \in \Real^{\infty \times p}$
is the quasimatrix whose $i$th function is $\KF(\cdot,x_i)$, and
$\vw \in \Real^p$ is a weight vector to be learned.
Further, $\nrm{\ff}_{\HK}^2 = \nrm{\vw}_{\K}^2$ (see \cref{sec:kernels-rkhs}).
Hence, defining
$\K \equiv \KF(\vx,\vx) = \QK(\vx) \in \Real^{p \times p}$,
\cref{eq:KRR} reduces to
\begin{equation}
  \label{eq:KRR-finite}
  \min_{\vw \in \Real^p}
  \fnrm{\vy - \K\vw}^2
  + \lambda \nrm{\vw}_{\K}^2,
\end{equation}
and the optimal solution is \cref{eq:KRR-solution-form} with weights $\vw$ given by
\begin{equation}
  \label{eq:KRR-solution}
  \vw = \prn{ \K + \lambda \Mx{I}{p}}^{-1} \vy.
\end{equation}

Kernel ridge regression
is a form of \emph{nonparametric regression} since we have not
assumed an underlying distribution for the function;
the user need only specify a kernel function which in turn specifies
the space of functions over which we will optimize.
Our solutions come from an infinite dimensional space $\HK$, yet the
minimization problem reduces to a simple least squares problem.
Full descriptions of RKHS properties can be found in many places in
the literature; see Wainwright \cite[Chapter 12]{Wa19} for one such
pedagogical treatment.
The remainder of this section provides a deeper
explanation of the properties of RKHS that make it amenable to
our problem.

\subsection{Background}
\label{sec:background}

Before proceeding, we recall some concepts important to our subsequent discussion.
The notation $L^2(\mathcal{I})$ denotes the Hilbert space of square integrable functions
on the interval $\mathcal{I}$ with inner product
$\langle f, g \rangle_{L^2(\mathcal{I})} \equiv \int_\mathcal{I} f(x) g(x) dx$.
We use $\ell^2(\mathbb{N})$ to denote the Hilbert
space of square-summable real-valued sequences,
i.e., $\ell^2(\mathbb{N}) \equiv \crly*{ \set{\theta_i}_{i\in\Natural} |
  \sum_{i=1}^{\infty} \theta_i^2  < \infty }$,
with inner product
$\langle \theta, \gamma \rangle_{\ell^2(\mathbb{N})} \equiv \sum_{i=1}^{\infty} \theta_i \gamma_i$.

\subsection{Kernels and RKHS} 
\label{sec:kernels-rkhs}
As mentioned above, one way to define an RKHS is to start with a
symmetric psd kernel function
The RKHS $\HK$ is then defined as
the closure of the set of linear combinations of
functions of the kernel, i.e.,
\begin{equation}\label{eq:HK-defn-1}
  \HK \equiv
  \overline{
  \crly*** {
    f = \sum_{i=1}^{n} w_i \KF(\cdot, x_i)
    \text{ for any } n \in \mathbb{N},
    \text{ points } \{x_i\}_{i=1}^n
    \text{, and } \Vc{w} \in \mathbb{R}^n
  } }.
\end{equation}
The intuition here is that the RKHS
contains all finite linear combinations of the univariate functions defined
by the kernel function at any set of design points.

For functions $f = \sum_{i=1}^{n} w_i \KF(\cdot, x_i) \in \HK$ and
$\tilde{f} =\sum_{j=1}^{\tilde{n}} \tilde{w}_j \KF(\cdot, \tilde{x}_j) \in \HK$,
the inner product is defined to be
\begin{equation}\label{eq:HK-innerprod-1}
  \ang*{ f, \tilde{f} }_{\HK} \equiv
  \sum_{i=1}^n \sum_{i=1}^{\tilde{n}} w_i \tilde{w}_j \KF(x_i, \tilde{x}_i).
\end{equation}
Hence, the norm in the regularization term of the kernel ridge regression problem
\cref{eq:KRR} takes the form:
\begin{equation*}
  \|f\|_{\HK}
  = \langle f, f \rangle_{\HK}
  = \sum_{i=1}^p \sum_{j=1}^{p} w_i w_j \KF(x_i, x_j)
  = \Vc{w}' \K \Vc{w} \equiv \|\Vc{w}\|_{\K}.
\end{equation*}
Using this inner product, $\HK$ has the \emph{kernel reproducing property}:
\begin{displaymath}
  \ang*{ f, \KF(\cdot, x) }_{\HK} = \sum_{i=1}^n w_i \KF(x_i, x) = f(x).
\end{displaymath}
We have not described the inner product and reproducing property for the elements
of $\HK$ obtained by taking the closure with respect to Cauchy sequences
and refer the interested reader to \cite[Theorem 12.11]{Wa19} for further details.

\subsection{Representer Theorem and Kernel Ridge Regression}
\label{sec:repr-theor-kern}

The kernel ridge regression problem in \cref{eq:KRR}
is tractable because of a special property of an RKHS.

\begin{theorem}[Representer Theorem]
  For a positive semi-definite kernel $\KF$, its associated RKHS
  $\HK$, and a set of $n$ observations $\set{(x_i, y_i)}_{i=1}^p$,
  consider the regularized empirical risk minimization problem defined
  by the loss function $\Qm{L}$:
  \begin{equation}
    \min_{\ff \in \HK} \sum_{i=1}^p \Qm{L} \prn*(y_i, \ff(x_i) ) + \lambda \|\ff\|_{\HK},
  \end{equation}
  where $\lambda > 0$.  The minimizer of this problem $\ff^*$ admits the form:
  \begin{equation}
    \ff^*= \sum_{j=1}^{p} w_j \KF(\cdot, x_j) \text{ for some } \Vc{w} \in \mathbb{R}^p .
  \end{equation}  
\end{theorem}

The remarkable aspect of the Representer Theorem is that it transforms
an infinite dimensional problem into a problem with only $p$ parameters.
In essence, the optimal
solution to the problem is a linear combination of just the $p$
univariate functions defined by the design points
$\{\KF(\cdot, x_1), \ldots, \KF(\cdot, x_d)\}$ while the RKHS as a
whole contains all finite linear combinations of the univariate
functions $\KF(\cdot, x)$ defined for any $x \in \mathcal{I}$.
Proof of these results and extended discussion can be found in \cite[Section 12.5.2]{Wa19}.

\subsection{Mercer's Theorem}
We next consider another useful characterization of the RKHS associated with a psd kernel in terms of the kernel's eigenfunction expansion (essentially an infinite dimensional feature map).
This is given by Mercer's theorem which
state that all positive semidefinite kernels admit such a
decomposition provided the kernel function is square integrable over
its domain, i.e.,
\begin{displaymath}
  \int_{\mathcal{I}} \int_{\mathcal{I}} \KF(x,y)^2 dx dy < \infty.
\end{displaymath}
This is also known as the Hilbert-Schmidt condition.

\begin{theorem}[Mercer's Theorem]
\label{thm:mercer}
  Let $\KF: \mathcal{I} \times \mathcal{I} \rightarrow \Real$ be a symmetric psd kernel.
  If $\KF$ is square integrable over its domain, then
  the kernel has an eigen-expansion of the form:
  \begin{equation}\label{eq:kernel-eigenexp}
    \KF(x, y) = \sum_{j=1}^{\infty} \mu_j \psi_j(x) \psi_j(y),
  \end{equation}
  where $\set{\psi_j(\cdot)}_{j=1}^{\infty}$ is a sequence of
  eigenfunctions that form an orthonormal basis for $L^2(\mathcal{I})$ and
  $\mu_1 \geq \mu_2 \geq \cdots \geq 0$ is a sequence of associated
  nonnegative eigenvalues.
  The infinite series in this eigen-expansion converges absolutely and uniformly.
  The eigenvalues and eigenfunctions satisfy
  \begin{displaymath}
    \int_{\mathcal{I}} \KF(x, y) \psi_j(y) dy = \mu_j \psi(x) \qtext{for all $j=1,2,\ldots$}
  \end{displaymath}
\end{theorem}

An extensive list of kernels and their associated Mercer series
representation can be found in Appendix A of \cite{FaMc15}.

Given such a set of eigenfunctions, an alternative 
characterization of $\HK$ and its associated inner product
can be given in terms of linear
combinations of the basis functions $\psi_j(x)$. 
Generally, one would not actually compute an eigenexpansion, but it
provides insight into the structure of the RKHS.
This can be contrasted with the definitions in \cref{eq:HK-defn-1,eq:HK-innerprod-1}.

\begin{corollary}
  Let $\KF: \mathcal{I} \times \mathcal{I} \rightarrow \Real$ be a symmetric psd kernel
  that is integrable over its domain.
  Using the eigen-expansion \cref{eq:kernel-eigenexp} from Mercer's theorem, 
  we have
  \begin{equation}\label{eq:HK-defn-2}
    \HK \equiv \crly***{ f = \sum_{{j=1}\atop {\mu_j \neq 0}}^{\infty} \beta_j \psi_j
      \; \bigg | \;
      \{\beta_j\}_{j=1}^{\infty} \in \ell^2(\mathbb{N})
      \text{ with } \sum_{{j=1}\atop{\mu_j \neq 0}}^{\infty} \frac{\beta_j^2}{\mu_j} < \infty
    },
  \end{equation}
  and the associated inner product for functions
  $f = \sum_{j=1}^{\infty} \beta_j \psi_j$ and $g = \sum_{j=1}^{\infty} \theta_j \psi_j$
  is
  \begin{equation}\label{eq:HK-innerprod-2}
    \langle f, g \rangle_{\mathcal{H}}
    \equiv \sum_{j=1}^{\infty} \frac{\langle f, \psi_j \rangle \langle g, \psi_j \rangle}{\mu_j}
    = \sum_{j=1}^{\infty} \frac{\beta_j \theta_j}{\mu_j}.
  \end{equation}
\end{corollary}

The RKHS norm used in the kernel ridge regression problem is
$\|f\|_{\HK}^2 = \langle f, f \rangle = \sum_{j=1}^{\infty}
\beta_j^2/\mu_j$.
Another perspective on the RKHS defined by a
Mercer kernel is that it is isomorphic to an elliptical subset of the
sequence space $\ell^2(\mathbb{N})$ defined by the non-negative
eigenvalues $\{\mu_j\}_{j=1}^{\infty}$.

\subsection{Broader Context}
A useful RKHS property is there exists ``universal kernels" for which the associated space can approximate continuous function arbitrarily well.
Both the exponential and Gaussian kernels are universal
\cite[Cor. 4.58]{StCh08}.
\begin{definition}[Universal Kernel \cite{MiXuZh06}]
  A continuous kernel function $\KF$ on $\mathcal{I}$ is said to be
  \emph{universal} if the associated RKHS $\HK$ is dense in $C(\mathcal{I})$,
  i.e. for every $g \in C(\mathcal{I})$ and any $\epsilon > 0$ there
  exists an $f \in \HK$ such that $\|f-g\|_{\infty} \leq \epsilon$.
\end{definition}

In our discussion, we have assumed the kernel function $\KF$ is directly specified; here we quickly survey several ways of constructing $\KF$ that result in several special cases of an RKHS.
The first is by instead specifying a feature map $\Phi(x): \Real \rightarrow \Real^D$ and constructing the kernel as $\KF(x, y) = \Phi(x)^{\Tr}\Phi(x)$.
One can confirm this will always produce a psd kernel function, and working with the
kernel can produce computational savings if the kernel can be computed in closed form (known colloquially as the ``kernel trick").
Thus the RKHS machinery can be used with any finite feature map or set of basis functions specified by the user.
Furthermore, this category also includes cubic splines and smoothing splines which can be represented by a finite basis set.

Another method is to essentially reverse Mercer's theorem to construct a kernel from a desired orthonormal basis for $L^2(\mathcal{I})$. Here the challenge is to select $\mu_j$ such that the resulting kernel satisfies the Hilbert-Schmidt condition and can be expressed in closed form.
For an example, see \cite[Section 3.9.2]{FaMc15} which gives two choices of $\{\mu_i\}_{i=1}^\infty$ to construct a kernel from the Chebyshev polynomials on $[-1,1]$.

Kernel ridge regression is closely connected to \emph{Gaussian process
regression}. A Gaussian process is a distribution over univariate
functions and is determined by a mean function
$\mu(x): \mathcal{I} \rightarrow \Real$ and a psd kernel function
$\KF: \mathcal{I} \times \mathcal{I} \rightarrow \Real$ that serves as the
covariance function. The mean function of this posterior
will be exactly the function outputted by kernel ridge regression; the
posterior covariance can also be given in closed form.

\section{Simple 3-Way Setting with One Continuous Mode}
\label{sec:simple-3-way}

To mitigate the overwhelming notation in the general $d$-way with an arbitrary
combination of discrete and continuous modes, we first establish the basic principles
using three-way example from \cref{eq:hcp-3way} with a single continuous mode; i.e.,
we have a quasitensor $\QT$ defined on $\mdom{m,n} \otimes \mathcal{I}$ where
$\mathcal{I}$ is a continuous interval in $\Real$.
We want to find a \cphifi decomposition of the form
$\QT \approx \dsqr{\A,\B,\QC}$ where $\A \in \Rmsiz{m,r}$,
$\B\in \Rmsiz{n,r}$, and $\QC \in \Rmsiz{\infty,r}$. This means
\begin{equation}\label{eq:hcp-3way-reduc}
  \QT(i,j,x) \approx \sum_{\ell=1}^r \va_{\ell}(i)\, \vb_{\ell}(j) \, \fc_\ell(x)
  \qtext{for all}
  (i,j,x) \in \mdom{m,n} \otimes \mathcal{I},
\end{equation}
where, for each $\ell \in [r]$, $\va{\ell}$ is the $\ell$th column of $\A$,
$\vb{\ell}$ is the $\ell$th column of $\B$,
and $\Fn{c}{\ell}$ is the $\ell$th function in the quasimatrix $\QC$.

To fit the model, we assume we have a discrete set of observations. 
We let $p$ denote the number of distinct values of $x$ such that  $\QT(i,j,x)$ is known for
at least one $(i,j) \in \mdom{m,n}$ and arrange these $p$ \emph{sample points} in mode 3 as
a vector
$\vv =
\begin{bmatrix}
  v_1 & v_2 & \cdots & v_p
\end{bmatrix}^{\Tr}
\in \Real^p$.
We define
\begin{equation}\label{eq:omega-3-way}
  \Omega = \set{ (i,j,k) \in \mdom{m,n,p} | \QT(i,j,v_k) \text{ is known}}.
\end{equation}
The optimal model will minimize the error between the model and the quasitensor at the observed points:
\begin{equation}\label{eq:opt-prob-3-v0}
  \min_{\A,\B,\QC} \frac{1}{2} \sum_{(i,j,k) \in \Omega} \prn**{ \QT(i,j,v_k) - \va_{\ell}(i)\, \vb_{\ell}(j) \, \fc_\ell(v_k) }^2.
\end{equation}

This can be expressed in more standard notation as follows.
Define the (discrete) sample tensor $\Tobs \in \Rmnp$ as
\begin{equation}\label{eq:Tobs-3-way}
  \Tobs(i,j,k) \equiv
  \begin{cases}
    \QT(i,j,\spt{k}) & \text{if } (i,j,k) \in \Omega \\
    0 & \text{if } (i,j,k) \not\in \Omega.    
  \end{cases}
\end{equation}
The tensor $\Tobs \in \Rmsiz{m,n,p}$ is the evaluation of $\QT\in\Rmsiz{m,n,\infty}$ at the vector $\vv \in \Real^p$ for the indices defined in $\Omega$.

We can write \cref{eq:opt-prob-3-v0} as
\begin{equation}\label{eq:opt-prob-3}
  \min_{\A,\B,\QC} \frac{1}{2} \nrm*{ \Tobs - \dsqr{\A,\B,\C} }_{\Omega}^2
  \qtext{subject to} \C = \QC(\vv).
\end{equation}
The matrix $\C \in \Rmsiz{p,r}$ is the evaluation of the quasimatrix
$\C \in \Rmsiz{\infty,r}$ at the vector $\vv \in \Rmsiz{p}$.

We consider two cases: (1) complete data with $\Omega = \mdom{m,n,p}$ (Fig.~\ref{fig:data}),
and (2) incomplete or misaligned data so that $\Omega \subset \mdom{m,n,p}$ (Fig.~\ref{fig:sampled-data}).
In the misaligned case, it may be that the sample points are distinct
for every $(i,j) \in \mdom{m,n}$.

To make determination of $\QC$ feasible,
we constrain the functions $\Fn{c}{\ell}$ to be from an RKHS
defined by positive semidefinite kernel function
$\KF: \mathcal{I} \times \mathcal{I} \rightarrow \Real$.
We assume the
\emph{design points} for the kernel are the same as the sample points,
and the design points are used to determine the set of functions that
will be combined linearly to form each function in $\QC$.
In other words, the quasimatrix $\KF(\cdot, \spt)$ is the dictionary
of functions that will be used to build up the functions in $\QC$ 
(because of the Representer Theorem; see \cref{sec:repr-theor-kern}).
We let $\W \in \Real^{p \times r}$ be the weight matrix such that
its $\ell$th column, $\Vc{w}{\ell}$, holds the weights for $\Fn{c}{\ell}$ so that
\begin{equation}\label{eq:basis-functions}
  \Fn{c}{\ell}(\cdot)
  = \msum{k/p} \Vc{w}{\ell}(k) \; \KF(\cdot, \spt)
  = \Qm[\hat]{K}^{\Tr} \Vc{w}{\ell} 
\end{equation}
where $\Qm[\hat]{K} = \KF(\cdot,\vv) \in \Real^{\infty \times p}$
denotes the quasimatrix such that its $k$th function is $\KF(\cdot,\spt)$.
Then we can write \cref{eq:basis-functions} simply as $\QC=\Qm[\hat]{K} \W$, i.e.,
a quasimatrix of size $\infty \times p$ times an $p \times r$ matrix.

From the design/sample points, we define the sample kernel matrix
$\K \in \Rmsiz{p,p}$
such that its $\K(i,j)=\KF(\spt{i},\spt{j})$.
By definition of the RKHS, the matrix $\K$ is positive semidefinite.
Using RKHS, we then have
\begin{displaymath}
  \C = \QC(\vv) = \Qm[\hat]{K}(\vv) \W = \KF(\vv,\vv) \W = \K\W.
\end{displaymath}
Hence, the \cphifi optimization problem \cref{eq:opt-prob-3} becomes
\begin{equation}\label{eq:opt-3-way}
  \begin{aligned}
  \min_{\A,\B,\W} \text{~~}&
  \frac{1}{2}
  \nrm*{\, \Tobs - \dsqr{\A,\B,\K\W} \,}_{\Omega}^2
  +
  \nrm*{\W}_{\K}^2\\
  \text{s.t.~~~} & 
  \A \in \Rmsiz{m,r}, 
  \B \in \Rmsiz{m,r}, \qtext[\;]{and}
  \W \in \Rmsiz{p,r}.
  \end{aligned}
\end{equation}
Here $\nrm{\cdot}_{\Omega}$ and $\nrm{\cdot}_{\K}$ are defined
in \cref{eq:omega-nrm-def,eq:mat-nrm-K}, respectively.
The regularization term is needed for the RKHS constraint on $\QC$.
We propose to solve this problem via an alternating optimization approach:
\begin{algorithmic}[1]
  \Repeat
  \State Solve \cref{eq:opt-3-way} for $\A$ with $\B$ and $\W$ fixed
  \State Solve \cref{eq:opt-3-way} for $\B$ with $\A$ and $\W$ fixed
  \State Solve \cref{eq:opt-3-way} for $\W$ with $\A$ and $\B$ fixed
  \Until{converged}
\end{algorithmic}
We describe the solution of the subproblems in the cases of complete
and incomplete data below.

\subsection{Complete data}
\label{sec:complete-data}

By complete data, 
we mean $\Omega = \mdom{m,n,p}$.

\subsubsection{Solving for the first or second  factor matrix}
\label{sec:solving-first-factor}

If we are solving for $\A$ or $\B$ with the other factors fixed, the problem reduces
to the standard CP-ALS problem. In other words, letting $\C=\K\W$, we need to minimize
\begin{displaymath}
   \frac{1}{2} \nrm+{ \Tobs - \dsqr*{\A,\B,\C} }^2,
\end{displaymath}
and the solution in well known \cite{KoBa09}.
The optimal $\A$ with $\B$ and $\C$ fixed is 
\begin{displaymath}
  \A = \prn*{ \Tm[\bar]{T}{1}(\C \krp \B) } \;\backslash\; (\C'\C \had \B'\B)
\end{displaymath}
using MATLAB-like notation.
Analogously, the optimal $\B$ with $\A$ and $\C$ fixed is 
\begin{displaymath}
  \B = \prn*{\Tm[\bar]{T}{2}(\C \krp \A)} \;\backslash\; (\C'\C \had \A'\A)  .
\end{displaymath}
In either case, assuming $r \ll \min\set{m,n,p}$,
the dominant cost is $\bigO(mnpr)$ for the MTTKRP.
\subsubsection{Solving for the third factor quasimatrix}
\label{sec:solving-third-factor}

If we are solving for $\W$ with $\A$ and $\B$ fixed,
the problem becomes more interesting. 
Holding $\A$ and $\B$ fixed, it is convenient to rewrite
the objective \cref{eq:opt-3-way}
using \cref{eq:nrm-def,eq:ktensor-unfold-3,eq:krn-vec,eq:omega-nrm-def} as
\begin{equation}\notag
  \min_{\W}
  \frac{1}{2} \,
  \nrm*{\sqr*{ (\B \krp \A) \krn \K } \vec(\W)
    - \vec\prn*(\Tm[\bar]{T}{3}) }_2^2 +
  \frac{\lambda}{2} \vec\prn(\W)' (\Mx{I}{r} \krn \K) \vec(\W).
\end{equation}
From this, we can calculate the gradient with respect to $\vec(\W)$
using \cref{eq:mc-1,eq:mc-2}.
Setting that gradient equal to zero and using \cref{eq:krn-transpose,eq:krn-krn,eq:krp-had} yields
\begin{align*}
  \Vc{0}
  &= \sqr*{ (\B \krp \A) \krn \K }^{\Tr}
  \prn**{ \sqr*{ (\B \krp \A) \krn \K } \vec(\W)
    - \vec\prn*(\Tm[\bar]{T}{3}) }
  + \lambda \prn*{ \Mx{I}{r} \krn \K } \vec(\W)
  \\
  &= 
  \sqr*{ (\B'\B \had \A'\A) \krn \K^2 + \lambda (\Mx{I}{r} \krn \K) } \vec(\W)
  - \sqr*{ (\B \krp \A) \krn \K }^{\Tr} \vec\prn*(\Tm[\bar]{T}{3})
  \\
  &= \prn{ \Mx{I}{r} \krn \K }
  \sqr**{
    \sqr*{ (\B'\B \had \A'\A) \krn \K + \lambda \Mx{I}{pr} } \vec(\W)
  - \sqr*{ \prn(\B \krp \A)' \krn \Mx{I}{p} } \vec\prn*(\Tm[\bar]{T}{3})
  }
\end{align*}
which (leaving out $\Mx{I}{r} \krn \K$) is true if
\begin{displaymath}
  \prn*{ (\B'\B \had \A'\A) \krn \K + \lambda \Mx{I}{pr}} \vec(\W)
  = \vec\prn*(\Tm[\bar]{T}{3}(\B \krp \A))
\end{displaymath}
For large enough $\lambda$, the matrix $(\B'\B \had \A'\A) \krn \K + \lambda \Mx{I}{pr}$
is positive definite so a solution must exist.
Computing the right-hand side is an MTTKRP and costs $O(mnpr)$ work.
Solving the linear system is $\bigO((rp)^3)$ work and will likely be the dominant cost if $p$ is of the same order as $m$ and $n$. This is more expensive than a usual CP-ALS solve because of the combination of the kernel matrix and the regularization.

\subsection{Incomplete or misaligned data}
\label{sec:incomplete-data}
By incomplete data, we mean $\Omega \subset \mdom{m,n,p}$ (strict subset).
Let $q = \abs{\Omega}$ be the number of observed entries.
We emphasize that we are not treating the missing entries as zero in the
optimization problem below because the $\nrm{\cdot}_{\Omega}$ measure ignores
entries corresponding to $(i,j,k) \not\in \Omega$.

\subsubsection{Solving for the first or second discrete factor}
\label{sec:solving-first-factor-incomp}

We consider the problem of  solving for $\B$ with the other factors fixed and let $\C=\K\W$.
(The problem of solving for $\A$ is analogous.)
Since we can omit the second term of \cref{eq:opt-3-way}
because it does not involve $\B$, our problem
reduces to
\begin{displaymath}
  \min_{\B} \frac{1}{2} \nrm+{  \Tobs - \dsqr*{\A,\B,\C} }_{\Omega}^2.
\end{displaymath}%
Using the tensor perfect shuffle matrix $\Mx{P}{2}$ from \cref{eq:tensor-perfect-shuffle}, the selection matrix $\Mx{S}{\Omega}$ from \cref{eq:omega-nrm-def}, the tensor unfolding in mode 2 from \cref{eq:ktensor-unfold-2}, and the Kronecker identity from \cref{eq:krn-vec}, 
we can rewrite the problem as a least squares problem in terms of $\vec(\B)$:
\begin{displaymath}
  \min_{\B} \frac{1}{2}
  \nrm+{
    \Mx{S}{\Omega}' \Mx{P}{2}' \sqr*{ (\C \krp \A) \krn \Mx{I}{n} } \vec\prn{\B}
    -
    \Mx{S}{\Omega}' \vec\prn{ \Tobs }
  }_2^2.
\end{displaymath}
The matrix $\Mx{S}{\Omega}' \Mx{P}{2}' \sqr*{ (\C \krp \A) \krn \Mx{I}{n} } $
has size $q \times rn$. Assuming $q > rn$,
solving this linear least squares via QR costs $O(qr^2n^2)$.
The problem can be decomposed into $n$ subproblems to reduce the
total cost to $\bigO(qr^2)$, but we omit these details since this is generally
not the dominant cost.
(Note that in the case where all modes were discrete this would correspond to the censored alternating least squares approach discussed by \cite{Ho23}.)

\subsubsection{Solving for the third continuous factor}
\label{sec:solving-third-factor-incomp}

This is only somewhat different than the case for full data.
We use the tensor perfect shuffle matrix $\Mx{P}{3}$ from \cref{eq:tensor-perfect-shuffle}, the selection matrix $\Mx{S}{\Omega}$ from \cref{eq:omega-nrm-def}, the tensor unfolding in mode 3 from \cref{eq:ktensor-unfold-3}, the Kronecker identity from \cref{eq:krn-vec},
and the norm identity from \cref{eq:mat-nrm-K}.
We also define $\Mx{Z} \equiv \B \krp \A$ to
rewrite the problem \cref{eq:opt-3-way} as a least squares problem in terms of $\vec(\W)$:
\begin{equation}\notag
  \min_{\W}
  \frac{1}{2} \,
  \nrm*{\Mx{S}{\Omega}' \Mx{P}{3} \prn{ \Mx{Z} \krn \K } \vec(\W)
    - \Mx{S}{\Omega}' \vec\prn(\Tobs) }_2^2 +
  \frac{\lambda}{2} \vec\prn(\W)' (\Mx{I}{r} \krn \K) \vec(\W).
\end{equation}
From this, we calculate the gradient with respect to $\vec(\W)$
using \cref{eq:mc-1,eq:mc-2}
and set that equal to zero to arrive at
\begin{displaymath}
  \sqr*{ \prn(\Mx{Z} \krn \K)' \Mx{P}{3}" \Mx{S}{\Omega}" \Mx{S}{\Omega}' \Mx{P}{3}'
  \prn(\Mx{Z} \krn \K) + \lambda (\Mx{I}{r} \krn \K) } \vec(\W)
  = \prn(\Mx{Z} \krn \K)' \Mx{P}{3}" \vec(\Tobs).
\end{displaymath}
In the above, we use the fact that $\Mx{S}{\Omega}" \Mx{S}{\Omega}' \vec(\Tobs) = \vec(\Tobs)$ since $\Tobs$ has entries corresponding to $(i,j,k) \not \in \Omega$ zeroed out.
Factoring out ($\Mx{I}{r} \krn \Mx{K}$), we have
\begin{displaymath}
  \sqr*{
    \underbrace{\prn(\Mx{Z} \krn \Mx{I}{p})' \Mx{P}{3}" \Mx{S}{\Omega}"}_{\Mx{Y}{1}'}
    \underbrace{\Mx{S}{\Omega}' \Mx{P}{3}' \prn(\Mx{Z} \krn \K)}_{\Mx{Y}{2}}
    + \lambda \Mx{I}{rp} } \vec(\W)
  = \vec\prn*(\Tm[\bar]{T}{3} \Mx{Z}).
\end{displaymath}
Exploiting this structure, forming $\Mx{Y}{1}$ and $\Mx{Y}{2}$ costs $\bigO{qpr}$.
Define $\Mx{Y} = \Mx{Y}{1}'\Mx{Y}{2}"$, which costs $\bigO(qp^2r^2)$ to form.
Then, we find $\W$ by solving
\begin{displaymath}
  (\Mx{Y} + \lambda \Mx{I}{rp}) \vec(\W) = \vec(\Tm[\bar]{T}{3} \Z).
\end{displaymath}
The right-hand side is an MTTKRP and costs $\bigO(mnpr)$.
The cost to solve the system is $\bigO(p^3r^3)$.
Unfortunately, this problem does not decompose as in the discrete case,
and this step will generally be the main bottleneck in \cphifi with incomplete or misaligned data.

\section{General $d$-Way Setting with Arbitrary Continuous Modes}
\label{sec:general-d-way}

We now consider a $d$-way quasitensor $\QT$ of defined on
a combination of infinite and finite dimensional modes.
We think of this tensor as having \emph{continuous} (infinite) and \emph{discrete} (finite) modes.
Let $\mathcal{C}$  and $\mathcal{D}$ partition $[d]$ into  continuous and discrete modes.  Without loss of generality, we assume that the continuous modes are ordered first in the tensor and that there are $c$ such modes, i.e., $|\mathcal{C}| = c$ and $|\mathcal{D}| = d - c$.
We have a quasitensor $\QT$ defined on
$\mathcal{I}_{1} \otimes \cdots \otimes \mathcal{I}_{c} \otimes [n_{c+1}] \otimes \cdots \otimes [n_{d}]$ where each $\mathcal{I}_k$ for $k \in \mathcal{C}$ is an interval on $\Real$.
We can think of $\QT$ as a function on $c$ real-valued inputs producing an output which is a tensor of order $(d-c)$.

Our goal is to find a \cphifi decomposition of the form
\begin{displaymath}
  \QT \approx \dsqr{\QA{1},\dots,\QA{c},\A{c+1},\dots,\A{d}}
  \in \Real^{ \infty \times \cdots \times \infty \times n_{c+1} \times \cdots \times n_d},
\end{displaymath}
where
$\QA{k} \in \Rmsiz{\infty,r}$ for $k \in \mathcal{C}$
and $\A{k} \in \Rmsiz{n_k,r}$ for $k \in \mathcal{D}$.

To fit the model, we assume we have a discrete set of observations.
For each continuous mode $k \in \mathcal{C}$, 
we let $n_k$ denote the number of distinct values of $\xk$ such that some
$\QT(\xk[1], \dots, \xk[c], i_{c+1}, \dots, i_d)$ is observed and arranges
those points into the vector
$\vk = \begin{bmatrix} \sptk{1} & \sptk{2} & \cdots & \sptk{n_k}\end{bmatrix} \in \Real^{n_k}$.
For each finite mode $k \in \mathcal{D}$, we assume without loss of
generality that the set of unique values of $i_k$ in
$\mathcal{V}$ is exactly $[n_k]$.
(If not, we would reduce $n_k$ and renumber the indices to make it true.)

Now that the observations in the infinite modes are indexed,
we can define a set of discrete indices for the observed points as
\begin{equation}\label{eq:omega-dway}
  \Omega = \set{ (\miwc) | \QT(\sptk[1]{i_1},\dots,\sptk[c]{i_c},i_{c+1},\ldots,i_d) \text{ is known} }.
\end{equation}
Finally, we can define a discrete sample tensor
$\Tobs \in \Rmsiz$ such that
\begin{displaymath}
  \Tobs(\miwc) =
  \begin{cases}
    \QT(\sptk[1]{i_1},\dots,\sptk[c]{i_c},i_{c+1},\ldots,i_d) & \text{if } (\miwc) \in \Omega \\
    0 & \text{if } (\miwc) \not\in \Omega.
  \end{cases}
\end{displaymath}
We can think of $\Tobs$ as the evaluation of $\T$ at the inputs
$\vk[1] \in \Real^{n_1}$ thought $\vk[c] \in \Real^{n_c}$, at the
indices specified by $\Omega$.

The optimal model will minimize the error between the model and the
quasitensor at the observed points:
\begin{displaymath}
  \min_{\crly{\QA{k}}_{k \in \mathcal{C}} \atop {\crly{\A{k}}_{k \in \mathcal{D}}}}
  \frac{1}{2} \sum_{(\miwc) \in \Omega} \prn***{
    \Tobs(\miwc) - \sum_{j=1}^r
    \prod_{k \in \mathcal{C}} \QA{k}(\sptk{i_k},j)
    \prod_{k \in \mathcal{D}} \A{k}(i_k,j)
  }^2.
\end{displaymath}
We can rewrite this as
\begin{displaymath}
  \min_{\crly{\QA{k}}_{k \in \mathcal{C}} \atop {\crly{\A{k}}_{k \in \mathcal{D}}}}
  \frac{1}{2} \nrm*{ \Tobs - \dsqr{\miwc[\A]} }^2_{\Omega}
  \qtext{subject to}
  \A{k} = \QA{k}(\vk) \qtext[\;]{for all} k \in \mathcal{C}.
\end{displaymath}
Each matrix $\A{k} \in \Real^{n_k \times r}$ for $k \in \mathcal{C}$
is the evaluation of the quasimatrix $\QA{k} \in \Real^{\infty \times r}$ at
the vector $\vk \in \Real^{n_k}$.

As before, we consider two cases: (1) complete data with
$\Omega = [n_1] \otimes \cdots \otimes [n_d]$ and (2) incomplete or
misaligned data so that
$\Omega \subset [n_1] \otimes \cdots \otimes [n_d]$.  The key insight
we will develop for both cases is that fitting the \cphifi
Decomposition in an alternating fashion will only require us to
consider whether the \emph{current} mode is discrete or continuous.

For each continuous mode $k \in \mathcal{C}$, we constrain the
functions $\fakl$ to be from an RKHS defined by a positive
semidefinite kernel function
$\KF[k]: \mathcal{I}_k \times \mathcal{I}_k \rightarrow \Real$.  We
set the design points to be the same as the sample points and shared
for all $\fakl$ with $\ell \in [r]$.  We define the weight matrix
$\W{k}$ such that its $\ell$th column, $\vwkl$, holds the weights for
$\fakl$. Then the $\ell$th functional in $\QA{k}$ is
\begin{equation}\label{eq:f-k-j}
  \fakl(x) =
  \msum{j/n_{k}} \vwkl (j) \; \KF[k](x, \sptk{j})
  \qtext{for all} (k,\ell) \in \mathcal{C} \otimes [r].
\end{equation}
We can let $\Qm[\hat]{K}{k} \in \Real^{\infty \times n_k}$
denote the quasimatrix such that its $j$th function is $\KF[k](\cdot,\sptk{j})$.
Then we can write \cref{eq:basis-functions} simply as $\QA{k}=\Qm[\hat]{K}{k} \W{k}$.

To summarize, our model requires several things to be defined up front:
\begin{align*}
  r &\in \mathbb{N}
  && \text{\footnotesize(model rank)},\\
  \KF[k] &: \mathcal{S}_k \times \mathcal{S}_k \rightarrow \Real
  && \text{\footnotesize(psd kernel function for each  $k\in\mathcal{C}$)},\\
  n_k &\in \mathbb{N}
  && \text{\footnotesize(number of design points for each  $k \in \mathcal{C}$), and}\\
  \set{\sptk{1},\dots,\sptk{n_k}} &\subset \mathcal{S}_{k}
  && \text{\footnotesize(set of design points for each  $k\in\mathcal{C}$)}.
\end{align*}
For each $k \in \mathcal{C}$, we define $\K{k}$ to be the $n_k \times n_k$ sample
kernel matrix %
such that
\begin{displaymath}
  \Kk(i,j) = \KF[k](\sptk{i},\sptk{j}) \qtext{for all}
  k \in \mathcal{C} \qtext{and} (i,j) \in \mdom{n_k,n_k}.
\end{displaymath}
Each sample kernel matrix $\K{k}$ for $k \in \mathcal{C}$ is then
defined by $\KF[k]$ and the design points
$\set{\sptk{1},\dots,\sptk{n_k}}$; these matrices are precomputed and
do not change while fitting the model.

The parts to be learned are
\begin{align*}
    \W{k} = \begin{bmatrix} \vwkl{1} & \vwkl{2} & \cdots & \vwkl{r} \end{bmatrix}
  &\in \Rmsiz{n_k,r}\qtext{for all} k \in \mathcal{C}
  \qtext{and}\\
  \A{k} = \begin{bmatrix} \vakl{1} & \vakl{2} & \cdots & \vakl{r} \end{bmatrix}
  &\in \Rmsiz{n_k,r}\qtext{for all} k \in \mathcal{D}
  .
\end{align*}
The $d$-way \cphifi optimization problem assumes we have
$\Tobs \in \Rmsiz$ and semipositive definite kernel matrices $\Kk$ of
size $n_k \times n_k$ for all $k \in \mathcal{C}$ and then is given
by:
\begin{equation}\label{eq:opt-d-way}
  \begin{aligned}    
  \min_{\set{\A{k}}_{k \in \mathcal{D}} \atop \set{\W{k}}_{k \in \mathcal{C}}}\text{~~~~} &
  \frac{1}{2}
  \nrm*{ \Tobs - \dsqr{\K{1} \W{1},\dots,\K{c} \W{c},\A{c+1},\dots,\A{d}} }^2_{\Omega}
  + \frac{\lambda}{2} \sum_{k \in \mathcal{C}} \nrm*{\W{k}}^2_{\K{k}} \\
  \text{subject to~~~~}
  &  \Mx{W}{k} \in \Rmsiz{n_k,r}
  \qtext[\;]{for all} k \in \mathcal{C}
  \qtext{and} \A{k} \in \Rmsiz{n_k,r} \qtext[\;]{for all} k \in \mathcal{D}.
  \end{aligned}
\end{equation}
We propose to solve this problem via an alternating optimization
approach in which the learned components for all modes except $k$ are
held fixed, and then we solve for $\A{k}$ if $k \in \mathcal{D}$ or
$\W{k}$ if $k \in \mathcal{C}$.  We then cycle through the modes in
this fashion until a convergence criterion is achieved.  We use the
following definitions throughout.  Since the $\W{\ell}$ matrices are
fixed for $\ell \neq k$, we let $\A{\ell} = \Kk[\ell]\W{\ell}$ for all
$\ell \in \mathcal{C}$ with $\ell \neq k$.  We also define
\begin{align}
  \label{eq:Zk}
  \Mx{Z}{k} &= \A{d} \krp \cdots \krp \A{k+1} \krp \A{k-1} \krp \cdots \krp \A{1}
  , \text{ and}\\
  \label{eq:Vk}
  \V{k} &= \A{d}'\A{d}" \had \cdots \had \A{k+1}'\A{k+1}" \had \A{k-1}'\A{k-1}" \had \cdots
  \had \A{1}'\A{1}"
  .
\end{align}
In the subsequent sections, we describe the solution of each of these subproblems in the case of complete and incomplete data.

\subsection{Complete Data}
\label{sec:complete-d-way}
By complete data we mean that $\Omega = [n_1] \otimes \cdots \otimes [n_k]$.

\subsubsection{Solving for a discrete mode factor matrix}
\label{sec:solve-discrete-factor}
Consider solving for $\A{k}$ with $k \in \mathcal{D}$.
Then \cref{eq:opt-d-way} reduces to the standard CP-ALS problem
\begin{displaymath}
  \min_{\A{k}} \nrm*{ \Tobs - \dsqr{ \miwc[\A] }}^2,
\end{displaymath}
and the solution is well known \cite{KoBa09}. The
optimal $\A{k}$ with all other factors fixed is given by
\begin{equation}\label{eq:grad-Ak}
  \A{k} = \Tm[\bar]{T}{k}  \Z{k} \backslash \V{k} 
\end{equation}
using MATLAB-like notation and where $\Z{k}$ and $\V{k}$ are defined in
\cref{eq:Zk,eq:Vk}.

The algorithm is shown in \cref{alg:finite-mode}.
\Lineref{line:finite-mode-Mk} calculates the MTTKRP, whose
computational complexity is $\bigO(r \prod_k n_k)$ and dominates the
other costs.  \Lineref{line:finite-mode-Vk} calculates the gram
matrices and their Hadamard products, for a cost of
$\bigO( \sum_{\ell \neq k} n_{\ell} r^2 )$.
\Lineref{line:finite-mode-Ak} solves the linear system for a cost of
$\bigO(n_kr^2)$.

\begin{algorithm}[ht]
  \caption{Solving finite-dimensional mode (complete data)}
  \label{alg:finite-mode}
  \begin{algorithmic}[1]
    \Require $k \in \mathcal{D}$,
    complete tensor $\Tobs \in \Rmsiz!$, and
    factor matrices $\set{\A{\ell}}_{\ell \neq k}$
    where $\A{\ell} = \Mx{K}{\ell}\Mx{W}{\ell}$ for $\ell \in \mathcal{C}$ 
    \Ensure $\A{k} = \arg \min_{\A{k}} \nrm{ \Tobs - \dsqr[\miwc![\A]] }^2$
    \Function{SolveFiniteMode}{$k$, $\Tobs$, $\set{\A{\ell}}_{\ell \neq k}$}
    \State\label{line:finite-mode-Mk}$\Mx{Y} \gets \Mx[\bar]{T}{k} \prn*{\A{d} \krp \cdots \krp \A{k+1} \krp \A{k-1} \krp \cdots \krp \A{1}}$
    \Comment MTTKRP
    \State\label{line:finite-mode-Vk}$\Mx{V} \gets \A{d}'\A{d}" \had \cdots \had \A{k+1}'\A{k+1}" \had \A{k-1}'\A{k-1}" \had \cdots  \had \A{1}'\A{1}"$
    \State\label{line:finite-mode-Ak}$\Mx{A}{k} \gets \Mx{Y} \, \backslash \Mx{V}$
    \Comment solution to linear system $\Mx{A}{k}\Mx{V} = \Mx{M}$
    \State \Return $\A{k}$
    \EndFunction
  \end{algorithmic}
\end{algorithm}
\begin{algorithm}[ht]
  \caption{Solving infinite-dimensional mode (complete data)}
  \label{alg:infinite-mode}
  \begin{algorithmic}[1]
    \Require $k \in \mathcal{C}$,
    complete tensor $\Tobs \in \Rmsiz!$, 
    kernel matrix $\K{k}$, and
    factor matrices $\set{\A{\ell}}_{\ell \neq k}$
    where $\A{\ell} = \Mx{K}{\ell}\Mx{W}{\ell}$ for $\ell \in \mathcal{C}$ 
    \Ensure $\W{k} = \arg \min_{\W{k}} \nrm{ \Tobs - \dsqr[\A{1},\dots,\A{k-1},\K{k}\W{k},\A{k+1},\dots,\A{d}] }^2$
    \Function{SolveInfiniteMode}{$k$, $\Tobs$, $\K{k}$, $\set{\A{\ell}}_{\ell \neq k}$}
    \State\label{line:infinite-mode-MTTKRP}$\Vc{y} \gets \vec\prn*{ \Mx[\bar]{T}{k} \prn{ \A{d} \krp \cdots \krp \A{k+1} \krp \A{k-1} \krp \cdots \krp \A{1} }}$
    \Comment MTTKRP
    \State\label{line:infinite-mode-Vk}$\Mx{V} \gets \A{d}'\A{d}" \had \cdots \had \A{k+1}'\A{k+1}" \had \A{k-1}'\A{k-1}" \had \cdots  \had \A{1}'\A{1}"$
    \State\label{line:infinite-mode-M}$\Mx{M} \gets \V \krn \K{k} + \lambda \Mx{I}{rn_k}$
    \State\label{line:infinite-mode-wk}$\Vc{w} \gets  \Mx{M} \,/\, \Vc{y} $
    \Comment solution to linear system $\Mx{M} \Vc{w} = \Vc{y}$
    \State $\W{k} \gets \reshape(\Vc{w},n_k,r)$
    \State \Return $\W{k}$
    \EndFunction
  \end{algorithmic}
\end{algorithm}

\subsubsection{Solving for any continuous factor}
\label{sec:solve-continuous-factor}
Consider solving for $\W{k}$ with all other variables fixed.
As above, we ignore the regularization terms that are irrelevant to our purposes.
Then  \cref{eq:opt-d-way} can be rewritten
using \cref{eq:krn-vec,eq:nrm-def,eq:ktensor-unfold,eq:mat-nrm-K},
as
\begin{equation}\label{eq:w-func}
  \min_{\W{k}}
  \frac{1}{2} \tnrm*{\, \sqr*[ \Z{k} \krn \K{k} ] \vec\prn*(\W{k})
    - \vec\prn*( \Tm[\bar]{T}{k}) \, }^2
  + \frac{\lambda}{2} \vec(\W{k})^{\Tr} (\I_r \krn \K{k}) \vec(\W{k}).
\end{equation}
Here $\Z{k}$ is the Khatri-Rao product in \cref{eq:Zk}.
Computing the gradient and setting it equal to zero yields
\begin{equation}\notag
  \sqr*{ \V{k} \krn \K{k}^2 + \lambda (\I_r \krn \K{k}) }\vec(\W{k})
  = \sqr*[ \Z{k}' \krn \K{k} ]\vec\prn*( \Tm[\bar]{T}{k}).
\end{equation}

We factor out $\I_r \krn \K{k}$ to get a linear system of size $rn_k \times rn_k$:
\begin{equation}\notag  
  \sqr*{ \V{k} \krn \K{k} + \lambda \I_{rn_k} }\vec(\W{k})
  = \vec\prn*( \Tm[\bar]{T}{k}  \Z{k}) .
\end{equation}
The filtering out is important to get a multiple of the identity for regularization.
Computing the right-hand size is an MTTKRP and costs $O(r \prod_k n_k)$ work.
Solving the linear system is $(rn_k)^3$ work.
This is more expensive than a usual CP-ALS solve because of the combination
of the kernel matrix and the regularization. 

The algorithm is shown in \cref{alg:infinite-mode}.
The two most expensive lines are as follows.
\Lineref{line:infinite-mode-MTTKRP} computes the MTTKRP for a cost of
$\bigO(r \prod_{k=1}^d n_k)$.
\Lineref{line:infinite-mode-wk} solves the linear system for a cost of $\bigO(r^3n_k^3)$.

\subsection{Incomplete Data}
\label{sec:incomplete-d-way}
By incomplete data we mean that $\Omega \subset [n_1] \otimes \cdots \otimes [n_k]$.
We define $q \equiv \abs{\Omega}$.

\subsubsection{Solving for any discrete factor}
\label{sec:solve-discrete-factor-incomplete}

We consider the problem of  solving for $\A{k}$ for some fixed $k \in \mathcal{D}$.
Since we can omit the regularization term, \cref{eq:opt-d-way} 
reduces to
\begin{displaymath}
  \min_{\A{k}} \frac{1}{2} \nrm+{ \Mx{S}{\Omega}' \vec\prn*{ \Tobs - \dsqr*{\miwc[\A]} } }_2^2.
\end{displaymath}
We can rewrite the problem as a least squares problem in terms of $\vec(\A{k})$, where $\Z{k}$ is defined in \cref{eq:Zk} and $\Mx{S}{\Omega}$ is again the operator that picks out the entries in $\Omega$:
\begin{displaymath}
  \min_{\A{k}} \frac{1}{2}
  \nrm+{
    \Mx{S}{\Omega}' \Mx{P}{k}' \sqr*{ \Z{k} \krn \Mx{I}{n_k} } \vec\prn{\A{k}}
    -
    \Mx{S}{\Omega}' \vec\prn{ \Tobs }
  }^2.
\end{displaymath}
The matrix $\Mx{M} = \Mx{S}{\Omega}' \Mx{P}{k}' \sqr*{ \Z{k} \krn \Mx{I}{n_k} } \in \Real^{q \times rn_k}$ in Line 3 should not be formed via explicit matrix multiplication.
Each row in $\Mx{M}$ corresponds to a $(\miwc) \in \Omega$
and can be formed directly as
\begin{displaymath}
  \sqr{ \A{d}(i_d,:) \had \cdots \had \A{k+1}(i_{k+1}) \had \A{k-1}(i_{k-1},:) \had \cdots \had \A{d}(i_d,:) } \krn \Mx{I}{n_k}(i_k,:)
\end{displaymath}
for a total of $\bigO(qrn_k)$ work.
The cost to solve this is $\bigO(qn_k^2r^2)$. 
This problem has special structure such that it can be decomposed
into one problem per row of $\A{k}$, lowering the total cost to $\bigO(qr^2)$.
However, we omit the details since this is not the dominant cost.

\subsubsection{Solving for any continuous factor}
\label{sec:solve-continuous-factor-incomplete}

This is only somewhat different than the case for full data.
We will consider solving for $\W{k}$ for some fixed $k \in \mathcal{C}$
while holding $\A{\ell}$ for $\ell \neq k$ (and thus $\Z{k}$) fixed.
We can rewrite \cref{eq:opt-d-way}
as
\begin{equation}\notag
  \min_{\W{k}}
  \frac{1}{2} \,
  \nrm*{\Mx{S}{\Omega}' \Mx{P}{k} \prn{ \Mx{Z}{k} \krn \K{k} } \vec(\W{k})
    - \Mx{S}{\Omega}' \vec\prn(\Tobs) }_2^2 +
  \frac{\lambda}{2} \vec\prn(\W{k})' (\Mx{I}{r} \krn \K{k}) \vec(\W{k}).
\end{equation}
From this, we calculate the gradient with respect to $\vec(\W{k})$
and set that equal to zero to arrive at
\begin{displaymath}
  \sqr*{ \prn(\Mx{Z}{k} \krn \K{k})' \Mx{P}{k}" \Mx{S}{\Omega}" \Mx{S}{\Omega}' \Mx{P}{k}'
  \prn(\Mx{Z}{k} \krn \K{k}) + \lambda (\Mx{I}{r} \krn \K{k}) } \vec(\W{k})
  = \prn(\Mx{Z}{k} \krn \K{k})' \Mx{P}{k}" \vec(\Tobs).
\end{displaymath}
In the above, we use the fact that $\Mx{S}{\Omega}" \Mx{S}{\Omega}' \vec(\Tobs) = \vec(\Tobs)$ because $\Tobs$ already has entries zeroed out for indices not in $\Omega$.
Factoring out $\Mx{I}{r} \krn \Mx{K}{k}$, we have the linear system
\begin{displaymath}
  \sqr*{
    \underbrace{\prn(\Mx{Z}{k} \krn \Mx{I}{n_k})' \Mx{P}{k}" \Mx{S}{\Omega}"}_{\Mx{Y}{1}'}
    \underbrace{\Mx{S}{\Omega}' \Mx{P}{k}' \prn(\Mx{Z}{k} \krn \K{k})}_{\Mx{Y}{2}}
    + \lambda \Mx{I}{rn_k} } \vec(\W{k})
    = \vec\prn*(\Tm[\bar]{T}{k} \Mx{Z}{k}).
\end{displaymath}

Some comments on complexity:
\begin{itemize}
\item The right-hand side is an MTTKRP and costs $\bigO(r \prod_{k' = 1}^d n_k')$.
\item To form $\Mx{Y}{2}$, we want to avoid computing the product
  $\Mx{Z}{k} \krn \Mx{K}$, which would cost $\bigO{r n_k^2 \prod_{k' \neq k} n_k}$ operations
  plus $\prod_{k' = 1}^d n_k \times n_k r$ storage.
  Instead, we form only the $q$ rows selected by $\Mx{S}{\Omega}$.
  The row of $\Mx{Y}{2}$ picked out by $(i_1, i_2, \ldots, i_d) \in \Omega$ is
  \begin{displaymath}
    \sqr{ \A{1}(i_1,:)\had \cdots \had \A{k-1}(i_{k-1},:) \had \A{i+1}(i_{k+1},:) \had \cdots \had \A{d}(i_d,:) } \krn \K(i_k,:).
  \end{displaymath}
  The total cost to form $\Mx{Y}{2}$ is $\bigO(qr n_k )$ plus $q \times r n_k $ storage.
\item Similarly, to form $\Mx{Y}{1}$, we want to avoid computing
  $\Mx{Z}{k} \krn \Mx{I}{n_k}$, which would cost $\bigO{r \prod_{k' \neq k} n_{k'}}$ operations
  plus $\left ( \prod_{k' =1}^d n_{k'} \right ) \times \left (rn_k \right )$ storage.
  Instead, we form only the $q$ rows selected by $\Mx{S}{\Omega}$.
  The row of $\Mx{Y}{1}$ picked out by $(i_1, i_2, \ldots, i_d) \in \Omega$ is
  \begin{displaymath}
    \sqr{ \A{1}(i_1,:)\had \cdots \had \A{k-1}(i_{k-1},:) \had \A{k+1}(i_{k+1},:) \had \cdots \had \A{d}(i_d,:) } \krn \I_{n_k}(i_k,:).
  \end{displaymath}
  The total cost to form $\Mx{Y}{1}$ is $\bigO(qr)$ plus $q \times r n_k$ storage.
\item The cost to compute $\Mx{Y}{1}'\Mx{Y}{2}" + \lambda \Mx{I}{rn_k}$ is $\bigO(n_k^2qr)$.
\item The cost to solve the system is $\bigO(n_k^3r^3)$.
  Unfortunately, this problem does not decompose as in the discrete case,
  and this step will generally be the main bottleneck in \cphifi with incomplete or misaligned data.
\end{itemize}

\section{Numerical Examples}
\label{sec:numerical-examples}

We compare CP and \cphifi (both alternating least squares)
on a rank-3 synthetic quasi-tensor of size
$\msiz{4,3,\infty}$ whose decomposition is shown in
\cref{fig:cp-data}.
We denote the quasi-tensor as $\QT$ and can view it as
$\QT: \set{1,2,3} \otimes \set{1,2,3,4} \otimes [0,1] \rightarrow \Real$.
Alternatively, we can view this as a $\msiz{3,4}$ set of functions such
that $\QT_{ij}: [0,1] \rightarrow \Real$ for each $(i,j) \in \mdom{m,n}$.

The experiments differ in the number and alignment
of the samples used to form the observation tensor $\Tobs$.
The samples are noisy, which we constructed by evaluating $\T$ at the
sample and adding independent Gaussian noise.
All the factors of the original quasitensor are nonnegative, so we constrain the solutions
to be nonnegative using a standard nonnegative least squares solver
(\texttt{lsqnonneq} in MATLAB).
We run each method five times and use the best solution
in terms of objective function on the observed data.

\paragraph{Experiments \#1 and \#2: Dense aligned samples in temporal dimension}
\label{sec:exp-2}

\begin{figure}[ht]
  \centering
  \includegraphics{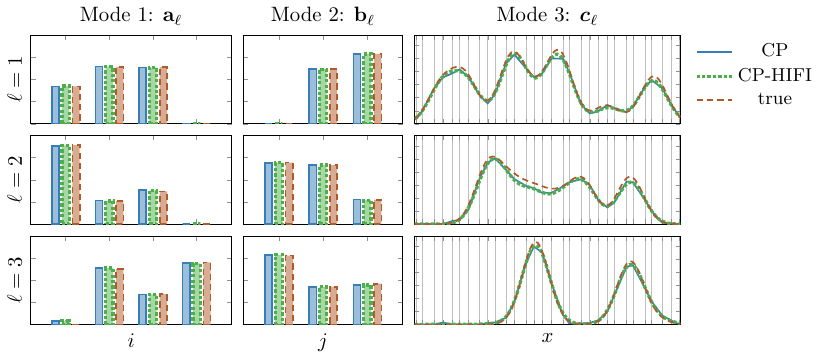}
  \caption{\textit{(Experiment 1)} 30 aligned points}
  \label{fig:exp-1}
\end{figure}

In the first experiment,
we evaluate each $\QT_{ij}$ at a set of
30 aligned $x$-values on the interval $[0,1]$, i.e., we evaluate the same
points from each of the 12 functions that constitute the quasitensor.
In this case, we have 360 data points arranged into the observation tensor
$\Tobs \in \Rmsiz{4,3,30}$ with no missing entries.
The resulting CP and \cphifi decompositions are shown in
\cref{fig:exp-1}.
The first column of plots are the columns of the factor matrix $\A \in \Rmsiz{4,3}$.
The second column of plots are the columns of the factor matrix $\B \in \Rmsiz{3,3}$.
The third column of plots is the columns of the factor matrix $\C$ for CP
and the functions from the factor quasi-matrix $\QC$ for \cphifi.
We have aligned the entries of $\C$ with the corresponding $x$-values
for the sample points.
The CP results are in blue with solid lines,
the \cphifi results are in green with dotted lines.
We also show the ``true'' underlying factorization that was used to
generate the data in brown dashed lines.
For  reference, the $x$-values corresponding to the sample points
are shown as vertical gray lines in the third column of plots.
The CP and \cphifi solutions are essentially the same
and very close to the true solution but inexact since the data are noisy.

In the second experiment, everything is the same except there are only 12 aligned
data points, so we have 144 data points arranged into the observation tensor
$\Tobs \in \Rmsiz{4,3,12}$.
The resulting CP and \cphifi decompositions are shown in
\cref{fig:exp-1}.
The CP and \cphifi solutions are essentially the same again, but the CP solution
is more jagged with fewer data points.

\begin{figure}[th]
  \centering
  \includegraphics{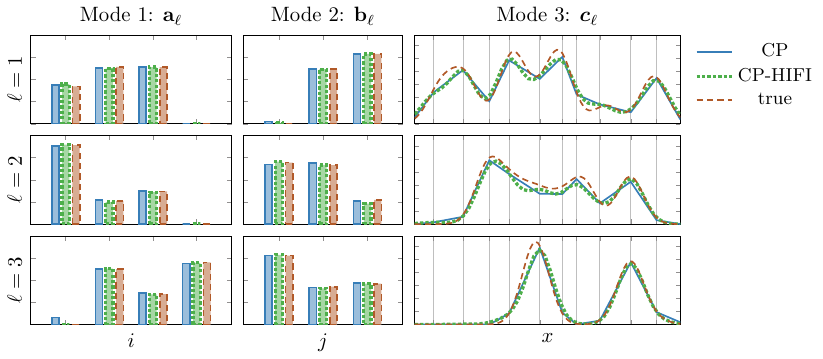}
  \caption{\textit{(Experiment 2)} 12 aligned points}
  \label{fig:exp-2}
\end{figure}

\paragraph{Experiment \#3: Dense aligned samples in temporal dimension with one unaligned point in one mode}
\label{sec:exp-3}

\begin{figure}[!th]
  \centering
  \includegraphics{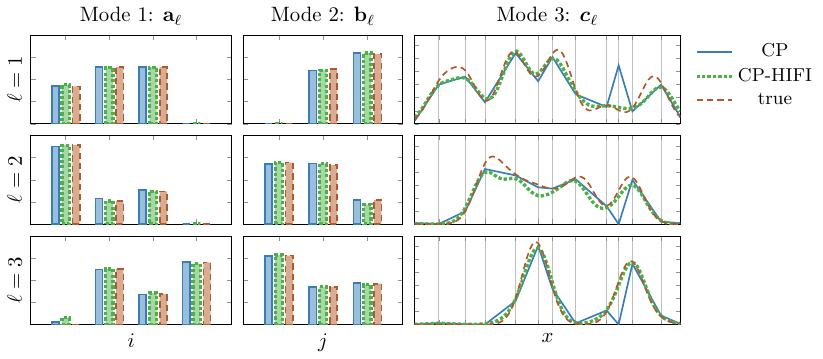}
  \caption{\textit{(Experiment 3)} 12 dense aligned samples for each $\QT_{ij}$ plus one extra
  data point for $\QT{22}$ at $x = 0.7677$}
  \label{fig:exp-3}
\end{figure}

In the third experiment, we evaluate each $\QT_{ij}$ at 12 aligned
points (same as experiment \# 2) on the interval $[0,1]$ \emph{plus} a
single extra data point for $\QT_{22}$ at $x = 0.7677$.
This means we have 145 data points in the data tensor
$\Tobs \in \Rmsiz{4,3,13}$ so that one slice
(corresponding to $x = 0.7677$)
has only a single known value and 11 missing values.
The results are shown in \cref{fig:exp-3}
Both the CP and \cphifi $\A$ and $\B$ factors are very close to the ground truth.
However, the $\C$ factor now has anomalous behavior
around $x = 0.7677$ since it has only a single observation
at that point.
Because the data is ``missing'' except for the single observed point
at $x=0.7677$, the resulting $\C$ matrix of CP shows spurious behavior
if we try to interpret it as a function.
We conclude that CP can be sensitive to even a single misaligned point.

\paragraph{Experiment \#4: Dense aligned samples in temporal dimension with one unaligned point in every mode}
\label{sec:exp-4}

\begin{figure}[th]
  \centering
  \includegraphics{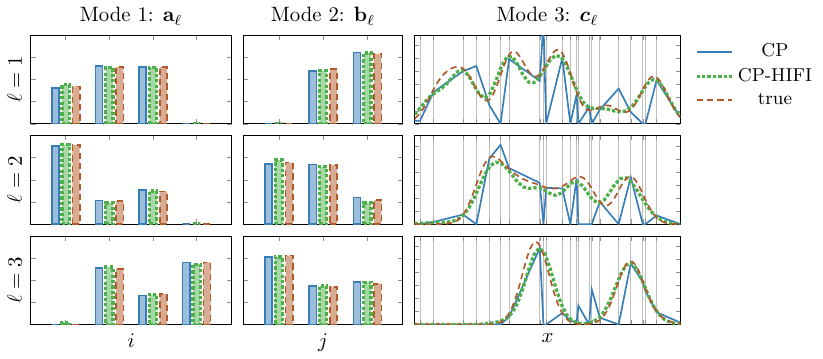}
  \caption{\textit{(Experiment 4)} 12 dense aligned samples for each $\QT_{ij}$ plus one extra
  data point for each $\QT{ij}$}
  \label{fig:exp-4}
\end{figure}

In the fourth experiment, we evaluate each $\QT_{ij}$ at 12 aligned
points (same as experiment \# 2) on the interval $[0,1]$ \emph{plus} a
single extra unaligned data point for each $\QT_{ij}$.
There are 156 data points.
The tensor $\Tobs$ is of size $\msiz{3,4,24}$, and there
are 12 slices that have only a single known entry.
The results are shown in \cref{fig:exp-4}.
The 24 different $x$-values where we have observed data are shown as gray bars
in the plots in the third column.
Both the CP and \cphifi $\A$ and $\B$ factors are very close to the ground truth.
However, the
resulting CP $\C$ factor is substantially different
from ground truth because of lots of spurious behavior.
More data degrades the solution.
In contrast, the \cphifi solution
is arguably improved as compared to experiment \#3
thanks to the additional data points.

\paragraph{Experiment \#5: Unaligned samples in temporal dimension}
\label{sec:exp-5}

\begin{figure}[th]
  \centering
  \includegraphics{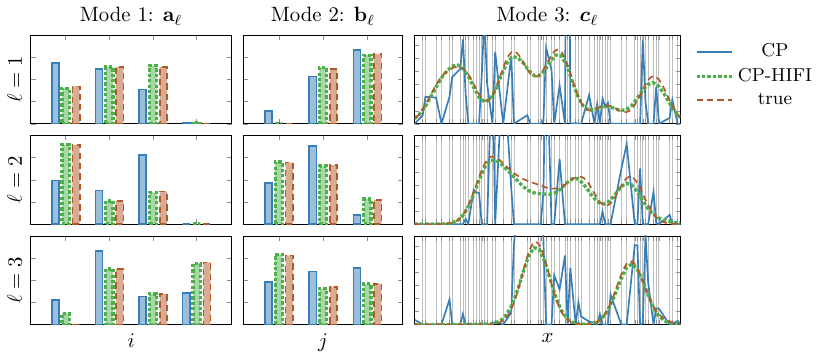}
  \caption{\textit{(Experiment 5)} 12 unaligned samples for each $\QT_{ij}$}
  \label{fig:exp-5}
\end{figure}

In the final experiment, we evaluate each $\QT_{ij}$ at 12 unaligned
points on the interval $[0,1]$ (\cref{fig:exp-5}).
This means we have a total of 144 observations (the same number as in experiment \#2),
but independently selected across fibers.
The resulting observation tensor $\Tobs$ is of size $\msiz{3,4,60}$, but with
only 2-3 known values per slice on average.
(We constrained the choice of $x$ values so that we
some repeats across the 144 observations; otherwise, there would be 144 slices
with a single known value per slice.)
The results are shown in \cref{fig:exp-5}.
The 60 different $x$-values are shown as gray lines in the third column of plots.
Now, the CP solution is significantly worse.
On the other hand, the \cphifi solution is arguably the
closest match we have seen. Rather than sowing chaos, the
diversity of observations have enables \cphifi to produce
a better fit to the
underlying smooth functions.

\section{Related work}
\label{sec:related-work}
Developing low-rank approximations which account for functional
structure has been the subject of interest across many fields of
mathematics and data science, including approximation theory, Bayesian
latent variable models, and stochastic processes.  Furthermore, many
applied fields such as neuroscience naturally have data in which one
or more of the modes is assumed to be observations from some smooth underlying function.
Here we review how methods developed across these fields relate to \cphifi.

\paragraph{Imposing Smoothness in CP}

In a standard CP model, constraints can be used to encourage
``smooth'' factors. In the discussion of related work that follows,  assume $\C \in \Rmsiz{p,r}$ is the factor matrix which should be smooth.

In \cite{TiKi02}, the authors propose constraining the smooth factor matrix
to be of the form $\C = \Mx{S}\Mx{W}$,
where $\Mx{S} \in \Rmsiz{p,p'}$ is a matrix of smooth b-splines and $\Mx{W} \in \Rmsiz{p',r}$ is a weight matrix to be learned.
In \cite[Chapter 4]{Ve18} and \cite{VeDeDe17}, the authors consider the same setup but for arbitrary known matrix $\Mx{S} \in \Rmsiz{p,p'}$ and discusses how $\Mx{S}$ can be a set of chosen basis functions.
Both are very similar to our constraint in terms of quasi-matrices $\QC = \Qm[\hat]{K}\W$
where $\Qm[\hat]{K}$ are the RKHS basis functions, but the representer theorem
guarantees the optimality of our solution in infinite dimensional space.

In \cite{MaSaVa08}, the authors propose a smoothness penalty of the form $\nrm{\Mx{L}{2}\C}^2$ where $\Mx{L}{2} \in \Rmsiz{(p-1),p}$ is the second difference matrix.
Others have adopted this approach such as \cite{Al13,YoZhCi16,KrLa23}.
(In \cite{KrLa23}, this is referred to as Whittaker smoothing.)

The framework discussed in this paper includes quasitensors for which all modes
are infinite dimensional, i.e., quasitensors which are a multivariate function defined on $\Real^3$:
$\QT(x,y,z)$. In this case, the goal of \cphifi is to 
approximate this as the sum of products of smooth functions:
\begin{equation}
  \Qt{T}(x,y,z) \approx \sum_{\ell=1}^r \fa_\ell(x) \, \fb_{\ell}(y) \, \fc_\ell(z)
  \qtext{for all}
  (x,y,z) \in \Real^3.
\end{equation}
Here, $\fa{\ell}$, $\fb{\ell}$, and $\fc{\ell}$ represent the $\ell$th
functions of quasimatrices $\QA,\QB,\QC \in \Real^{\infty \times r}$, respectively.
In the context of tensors, such a model was proposed
by Beylkin and Mohlenkamp \cite{BeMo02,BeMo05}.
The problem of fitting this decomposition is equivalent to functional regression in three dimensions with the assumption that the function factorizes across dimensions.

\paragraph{Continuous Analogues of (Multi-) Linear Algebra}

Various works consider analogues of matrix factorizations for quasimatrices \cite{TrBa97,Tr10,ToTr15}.
For example, Townsend and Trefethen \cite{ToTr15} define continuous analogues of the singular value decomposition (SVD) $\Qm{M} \approx \Qm{U}\Mx{\Sigma}\Mx{V}^*$ in which $\Qm{U}$ is  an $[a,b] \times r$ orthonormal quasimatrix, $\Mx{\Sigma}$ is an $r \times r$ diagonal matrix, and $\Mx{V}$ is an $r \times n$ orthonormal matrix such that
\begin{displaymath}
  \Qm{M}(x,j) \approx \sum_{\ell=1}^r \sigma_{\ell} \; \Fn{u}_{\ell}(x) \, \Vc{v}{\ell}(j).
\end{displaymath}
The functions $\set{\Fn{u}_{\ell}}$ in the quasimatrix $\Mx{U}$ are univariate functions on the interval $[a,b]$ represented by the coefficients of an expansion in Chebyshev polynomials (referred to as a \texttt{chebfun}). Orthogonality is defined with respect to the standard $L^2$ inner product.
Since we are restricted to the interval $[a,b]$, a truncated series can represent any continuous function to machine precision.
In addition, \cite{ShAv22} considers what they call semi-infinite linear regression (SILR) in which the design matrix is a quasimatrix; kernel ridge regression discussed in \Cref{sec:rkhs-repr-theor} is a special case of SILR.

The Karhunen–Lo\'eve decomposition, the infinite dimensional analogue of PCA, decomposes the covariance function of a quasimatrix $\Qm{X} \in \Real^{\infty \times n}$, i.e., $\Cm{C}(s,t) = \frac{1}{n} \sum_{i=1}^n (\Fn{x}{i}(t) - \Fn{\mu}(t) )(\Fn{x}{i}(s) - \Fn{\mu}(s))$ where $\Fn{\mu}(\cdot)$ is the mean function. As $\Cm{C}$ is a psd kernel function, this decomposition relies on Mercer's theorem. 

Townsend and Trefethen \cite{ToTr13,ToTr15} also consider a matrix with two continuous modes, i.e., a continuous bivariate function
$\Cm{M}: [a,b] \otimes [c,d] \rightarrow \mathbb{C}$,
which they term a \texttt{chebfun2}. They define QR, LU, and SVD for these objects.
In particular, the SVD can be used to find a rank-$r$ approximation that is accurate
to machine precision of the form
$\Cm{M} = \Qm{U}\Mx{\Sigma}\Qm{V}^*$ with $\Qm{U}$ as an $[a,b] \times r$ quasimatrix with orthonormal columns denotes $\Fn{u}{j}$, $\Mx{S}$ as an $r \times r$ diagonal matrix, and $\Qm{V}$ is an $[c,d] \times r$ quasimatrix with orthonormal columns denoted $\Fn{v}{j}$ such that
\begin{displaymath}
    \Cm{M}(x,y) \approx \sum_{j=1}^{\infty} \sigma_j \; \Fn{u}{j}(x) \, \Fn{v}{j}(y).
\end{displaymath}

Lastly, Hashemi and Trefethen \cite{HaTr17} consider trivariate functions over $[a,b] \times [c, d] \times [e,g]$ and introduce \texttt{chebfun3} to describe continuous third-order tensor and a continuous analogue of the Tucker decomposition. Let $\QT$ be the trivariate function,  $\Tn{G}$ be a discrete tensor of size $r_x \times r_y \times r_z$, $\Qm{U}$ be a quasimatrix of size $[a,b] \times r_x$, $\Qm{V}$ a quasimatrix of size $[c,d] \times r_y$, and $\Qm{W}$ a quasimatrix of size $[e,g] \times r_z$.  Then the \texttt{chebfun3} approximation is
\begin{displaymath}
    \QT(x,y,z) \approx \sum_{i=1}^{r_x} \sum_{j=1}^{r_y} \sum_{k=1}^{r_z} \Tn{G}(i,j,k) \; \Fn{u}{i}(x) \Fn{v}{j}(y) \Fn{w}{k}(z).
\end{displaymath}
As before, the component one-dimensional functions of the quasimatrices are represented by finite Chebyshev expansions.

\paragraph{Gaussian Process Factor Analysis (GPFA)}
Latent variable models are frequently used in analyzing complex, high-dimensional data where it is assumed there is a lower dimensional space of latent factors driving the variation in the data.
The Bayesian model is defined by a prior over the latent variables, and we fit the model to the observations by computing a posterior via expectation maximization (EM).
Consider the problem one observes a high-dimensional time series for which $m$ quantities are observed at $n$ time points and assembled into a  matrix $\Mx{Y} \in \Rmsiz{m, n}$ (e.g. neural spiking data collected from many neurons).
If we assume this high-dimensional data is actually being driven by a small number of low dimensional factors, we can incorporate low-rank structure and model our observed data as $\Mx{Y} \approx \Mx{W} \Mx{Z}$ where $r \ll m$,
Here $\Mx{Z} \in \Rmsiz{r,n}$ corresponds to a small number of latent trajectories
and $\Mx{W} \in \Rmsiz{m,r}$ is the ``loading" or ``readout" matrix relates the latents to the observed spikes.  
The observed time points $\{t_1, \ldots, t_n\}$ are assumed to be known and the same across all $n$ observations.

Probabilistic PCA, factor analysis, and Gaussian process Factor Analysis (GPFA) \cite{YuCuSaRy09} are all latent variable models that utilize this matrix factorization structure but differ in their choice of prior.
GPFA specifically imposes smoothness constraints over the latents $\Mx{Z}$ since the rows are assumed to be discrete observations of a continuous function over time. %
For a user-specified kernel $\KF$, GPFA uses the prior that the latents are drawn from the associated Gaussian process $z_j(t) \sim \mathcal{GP}(0, \KF)$ and then as in factor analysis allows a different variance for each neuron/observation: $y_{ij} \sim \mathcal{N}\big ([\Mx{W}\Mx{Z}]_{ij}, \sigma_i^2 \big )$.
In multi-output Gaussian process literature, this model has been called the intrinsic coregionalization model \cite[Chapter 4]{AlLoLq12}.

\section{Conclusions}
\label{sec:conclusions}
We propose \cphifi, a tensor decomposition tailored to the common real-world scenario of the data tensor being discrete measurements of a set of continuous processes.
The algorithm can accommodate any number of smooth modes, and we show how to efficiently work with misaligned measurements. 
Continuous modes are modeled using functions from an RKHS, which gives the practitioner 
the flexibility to work with a specific set of tailored basis functions or to use a 
universal kernel which can approximate any continuous function arbitrarily well.
Furthermore, the framework presented in this paper can easily be adapted to work with
other methods for functional regression.

As shown in the numerical experiments, even a single misaligned point can degredate
the factors found in the continuous mode. %
If there are only a few misaligned points, we can likely recover
the decomposition by adding an appropriate smoothness constraint. Such smoothing methods are
heuristic; whereas, the RKHS approach taken in this paper is both a computationally
efficient and principled way to fit functions to the continuous modes. %

The primary avenue we see for future work is the choice of design points. 
Currently the design points consist of all points in the continuous domain for which
at least one fiber has an observation. In cases with large-scale, misaligned tensors, 
it may be possible to effectively downsample the number of design points used for \cphifi.
In addition, it would be useful to give practitioners guidance on what points would be
most useful to sample in further data collection in order to improve the model.

\bibliographystyle{siamplainmod}


\end{document}